\documentclass[11pt]{amsart}
\usepackage{amsmath,amssymb,latexsym,amsmath,amsthm,amscd}

\theoremstyle{plain}
\newtheorem{thm}{Theorem}
\newtheorem{lem}[thm]{Lemma}
\newtheorem{prop}[thm]{Proposition}

\theoremstyle{definition}
\newtheorem{dfn}[thm]{Definition}
\newtheorem{ex}[thm]{Example}

\theoremstyle{remark}
\newtheorem{rmk}[thm]{Remark}

\DeclareMathOperator{\End}{End}

\DeclareMathOperator{\Ind}{Ind}

\begin{document}

\title{Vertex algebroids and Conformal vertex algebras associated with simple Leibniz algebras}
\author{Thuy Bui}
\address{Department of Management Science \& Information Systems \\ Rutgers Business School, Newark, NJ, USA}
\email{tb680@business.rutgers.edu}
\author{Gaywalee Yamskulna}
\address{Department of Mathematics\\ Illinois State University, Normal, IL, USA}
\email{gyamsku@ilstu.edu}


\begin{abstract} 
We first investigate the algebraic structure of vertex algebroids $B$ when $B$ are simple Leibniz algebras. Next, we use these vertex algebroids $B$ to construct indecomposable non-simple $C_2$-cofinite $\mathbb{N}$-graded vertex algebras $\overline{V_B}$. In addition, we classify $\mathbb{N}$-graded irreducible $\overline{V_B}$-modules and examine conformal vectors of these $\mathbb{N}$-graded vertex algebras $\overline{V_B}$.
\end{abstract}

\maketitle

\tableofcontents

\section{Introduction} 

\noindent The first aim of this paper is to classify vertex algebroids $B$ when $B$ are simple Leibniz algebras. The second aim of this paper is to construct irrational $C_2$-cofinite vertex algebras from the vertex algebroids. 

\vspace{0.2cm}

\noindent Vertex algebroids were introduced in 1999 in a series of studies on Gerbes of chiral differential operators in \cite{GMS} and on the chiral de Rham complex in \cite{MSV, MS1, MS2}. In these studies, Gorbounov, Malikov, Schechtman, and Vaintrob investigated the algebraic structure of $\mathbb{N}$-graded vertex algebras $V=\oplus_{n=0}^{\infty}V_{(n)}$  when $\dim V_{(0)}\geq 2$. They discovered that when they restricted the multiplication of $V$ to the vector space $V_{(0)}\oplus V_{(1)}$, they can produce many bilinear operations on $V_{(0)}\oplus V_{(1)}$ that give rise to many compatible relations. With these appropriate operations, one can show that $V_{(0)}$ is a unital commutative associative algebra, $V_{(1)} $ is a Leibniz algebra, $V_{(0)}\oplus V_{(1)}$ is a 1-truncated conformal algebra with additional relations. Later, all these compatible relations on $V_{(0)}\oplus V_{(1)}$ are summarized in the notion of vertex $V_{(0)}$-algebroid. It is worth to mention that for a given vertex $A$-algebroid $B$, one can construct an $\mathbb{N}$-graded vertex algebra $\oplus_{n=0}^{\infty}V_{(n)}$ such that $V_{(0)}=A$ and $V_{(1)}=B$ (cf. \cite{GMS}). Also, the classification of graded simple (twisted)-modules of vertex algebras associated with vertex algebroid were studied in \cite{LiY,LiY2} by Li and one of authors of this paper. Recently, Jitjankarn and one of authors of this paper investigated vertex algebroids associated with (semi)simple Leibniz algebras that have the simple Lie algebra $sl_2$ as their Levi factor (\cite{JY}), and constructed an indecomposable non-simple $C_2$-cofinite vertex algebra from a certain vertex algebroid such that its Levi factor is isomorphic to $sl_2$. In addition, they showed that this vertex algebra has only two irreducible $\mathbb{N}$ graded-modules (\cite{JY2}). Clearly, vertex algebroid plays an important role in the study of representation theory of $\mathbb{N}$-graded vertex algebras. However, the classification of vertex algebroids is far from being completed. Therefore, it is fundamentally important for us to examine vertex algebroids and their modules. 

\vspace{0.2cm}

\noindent The first part of this paper (Sections 2-3) is devoted to the study of the rich algebraic structure of vertex algebroid. We will first generalized the results in \cite{JY} by studying vertex $A$-algeroids $B$ as modules of the simple Lie algebra $sl_2$.  In particular, we show that if $sl_2$ is contained in $B$ then the unital commutative associative algebra $A$ is a local algebra. Next, we construct vertex algebroids $B_{\mathfrak{g}}$ from unital commutative associative algebras and Lie algebras $\mathfrak{g}$. We then examine the algebraic structure of the constructed vertex algebroid $B_{\mathfrak{g}}$ when $\mathfrak{g}$ are finite dimensional simple Lie algebras. Also, we analyze the algebraic structure of arbitrary vertex algebroids $B$ when $B$ are simple Leibniz algebras and prove that, under a suitable condition, the vertex algebroids $B$ are equivalent to vertex algebroids $B_{\mathfrak{g}}$ that are previously constructed. 

\vspace{0.2cm}

\noindent The $C_2$-cofiniteness property of vertex (operator) algebras plays a crucial role in the study of representation theory of vertex (operator) algebras 
(e.g. \cite{A2}, \cite{ ABD}, \cite{Bu}, \cite{DLM}, \cite{GN}, \cite{Mi1}-\cite{Mi3}, \cite{ Z}). Over the years, a large portion in the literature on the study of vertex (operator) algebras and their representations from a mathematical and physical point of view has been devoted to rational $C_2$-cofinite vertex algebras (e.g. \cite{B1}-\cite{B2}, \cite{D}, \cite{DoL}, \cite{DLM}, \cite{DoM1}, \cite{FLM2}, \cite{FZ}). However, the relatively recent discovery of aspects of logarithmic conformal field theory has changed the landscape of the representation theory of vertex (operator) algebras \cite{CR}. Mathematicians and physicists have begun to pay closer attention to irrational $C_2$-cofinite vertex algebras. However, there are few known examples of families of irrational $C_2$-cofinite vertex (super) algebras (e.g. \cite{A}, \cite{AdM1}-\cite{AdM4}, \cite{CF}, \cite{FFHST}, \cite{FGST1}-\cite{FGST3}, \cite{JY2}). Therefore, it is crucial for us to search for more families of irrational $C_2$-cofinite vertex algebras.

\vspace{0.2cm}

\noindent In the second part of this paper (Section 4),  we investigate a new family of irrational $C_2$-cofinite vertex algebras. Precisely, we construct indecomposable $C_2$-cofinite $\mathbb{N}$-graded vertex algebras from the vertex algebroids associate with simple Leibniz algebras that were obtained in the first part of this paper. In addition, we classify irreducible $\mathbb{N}$-graded modules of these vertex algebras and show a 1-1 correspondence between these $\mathbb{N}$-graded irreducible modules and irreducible modules of rational affine vertex operator algebras. We end this paper by examining conformal vectors of the indecomposable $C_2$-cofinite $\mathbb{N}$-graded vertex algebras associated with simple Leibniz algebras.

\section{A study of vertex algebroids as modules of simple Lie algebra $sl_2$} 

\noindent In subsection 2.1, we provide necessary background on vertex algebroids. In subsection 2.2, we study the algebraic structure of vertex algebroids when we consider these vertex algebroids as $sl_2$-modules.

\subsection{Preliminaries on 1-Truncated Conformal Algebras and Vertex Algebroids}

\begin{dfn}\cite{GMS} A {\em 1-truncated conformal algebra} is a graded vector space $C=C_0\oplus C_1$ equipped with a linear map $\partial:C_0\rightarrow C_1$ and bilinear operations $(u,v)\mapsto u_iv$ for $i=0,1$ of degree $-i-1$ on $C=C_0\oplus C_1$ such that the following axioms hold:

\medskip

\noindent(Derivation) for $a\in C_0$, $u\in C_1$,
\begin{equation}
(\partial a)_0=0,\ \ (\partial a)_1=-a_0,\ \ \partial(u_0a)=u_0\partial a;
\end{equation}

\noindent(Commutativity) for $a\in C_0$, $u,v\in C_1$,
\begin{equation} 
u_0a=-a_0u,\ \ u_0v=-v_0u+\partial(u_1v),\ \ u_1v=v_1u;
\end{equation}

\noindent(Associativity) for $\alpha,\beta,\gamma\in C_0\oplus C_1$,
\begin{equation}
\alpha_0\beta_i\gamma=\beta_i\alpha_0\gamma+(\alpha_0\beta)_i\gamma.
\end{equation}
\end{dfn}

\begin{dfn}(\cite{Br1}, \cite{Br2}, \cite{GMS}) Let $(A,*)$ be a unital commutative associative algebra over $\mathbb{C}$ with the identity $1$. A {\em vertex $A$-algebroid} is a $\mathbb{C}$-vector space $\Gamma$ equipped with 
\begin{enumerate}
\item a $\mathbb{C}$-bilinear map $A\times \Gamma\rightarrow \Gamma, \ \ (a,v)\mapsto a\cdot v$ such that $1\cdot v=v$ (i.e. a nonassociative unital $A$-module),
\item a structure of a Leibniz $\mathbb{C}$-algebra $[~,~]:\Gamma\times \Gamma\rightarrow\Gamma$, 
\item a homomorphism of Leibniz $\mathbb{C}$-algebra $\pi:\Gamma\rightarrow Der(A)$,
\item a symmetric $\mathbb{C}$-bilinear pairing $\langle ~,~\rangle:\Gamma\otimes_{\mathbb{C}}\Gamma\rightarrow A$,
\item a $\mathbb{C}$-linear map $\partial :A\rightarrow \Gamma$ such that $\pi\circ \partial =0$ which satisfying the following conditions:
\begin{eqnarray*}
&&a\cdot (a'\cdot v)-(a*a')\cdot v=\pi(v)(a)\cdot \partial(a')+\pi(v)(a')\cdot \partial(a),\\
&&[u,a\cdot v]=\pi(u)(a)\cdot v+a\cdot [u,v],\\
&&[u,v]+[v,u]=\partial(\langle u,v\rangle),\\
&&\pi(a\cdot v)=a\pi(v),\\
&&\langle a\cdot u,v\rangle=a*\langle u,v\rangle-\pi(u)(\pi(v)(a)),\\
&&\pi(v)(\langle v_1,v_2\rangle)=\langle [v,v_1],v_2\rangle+\langle v_1,[v,v_2]\rangle,\\
&&\partial(a*a')=a\cdot \partial(a')+a'\cdot\partial(a),\\
&&[v,\partial(a)]=\partial(\pi(v)(a)),\\
&&\langle v,\partial(a)\rangle=\pi(v)(a)
\end{eqnarray*}
for $a,a'\in A$, $u,v,v_1,v_2\in\Gamma$.
\end{enumerate}
\end{dfn}

\begin{prop}\cite{LiY} Let $(A,*)$ be a unital commutative associative algebra and let $B$ be a module for $A$ as a nonassociative algebra . Then a vertex $A$-algebroid structure on $B$ exactly amounts to a 1-truncated conformal algebra structure on $C=A\oplus B$ with 
\begin{eqnarray*}
&&a_ia'=0,\\
&&u_0v=[u,v],~u_1v=\langle u,v\rangle,\\
&&u_0a=\pi(u)(a),~ a_0u=-u_0a
\end{eqnarray*} for $a,a'\in A$, $u,v\in B$, $i=0,1$ such that 
\begin{eqnarray*}
&&a\cdot(a'\cdot u)-(a*a')\cdot u=(u_0a)\cdot \partial a'+(u_0a')\cdot \partial a,\\
&&u_0(a\cdot v)-a\cdot (u_0v)=(u_0a)\cdot v,\\
&&u_0(a*a')=a*(u_0a')+(u_0a)*a',\\
&&a_0(a'\cdot v)=a'*(a_0v),\\
&&(a\cdot u)_1v=a*(u_1v)-u_0v_0a,\\
&&\partial(a*a')=a\cdot \partial(a')+a'\cdot \partial(a).
\end{eqnarray*}
\end{prop}

\begin{prop}\cite{JY} Let $(A,*)$ be a finite dimensional commutative associative algebra and let $B$ be a finite dimensional vertex $A$-algebroid. Then 

\vspace{0.2cm}

\noindent (i) $Ker(\partial)=\{a\in A~|~\partial(a)=0\}$ is a subalgebra of $A$ and $B$ acts trivially on $Ker(\partial)$.

\vspace{0.2cm}

\noindent (ii) $Ker(\partial)$ contains every idempotent of $A$.

\vspace{0.2cm}

\noindent (iii) $\partial(A)$ is a $Ker(\partial)$-module.
\end{prop}

 
\subsection{Vertex algebroids as $sl_2$-modules}

\ \ 

\vspace{0.2cm}

\noindent Let $(A,*)$ be a finite dimensional commutative associative algebra with the identity $\hat{1}$ such that $\dim A\geq 2$. Let $B$ be a vertex $A$-algebroid. Assume that there exist $e,f,h\in B$ such that $e_0f=h$, $h_0e=2e$, $h_0f=-2f$, $h_0h=0$ and $Span\{e,f,h\}$ is a Lie algebra that is isomorphic to $sl_2$. In this subsection, we will investigate the algebraic structures of the commutative associative $A$ and the vertex $A$-algebroid $B$ as $sl_2$-modules. Also, under suitable conditions, we will show that $A$ is a local algebra. 

\vspace{0.1cm}

\noindent Because $sl_2$ is semisimple and $\mathbb{C}{\hat{1}}$ is a $sl_2$-submodule of $A$, there exists an $sl_2$-module $N$ such that 
\begin{eqnarray*}
A=\mathbb{C}{\hat{1}}\oplus N.
\end{eqnarray*} 
We set $N\cdot B=Span\{n\cdot b~|~n\in N, b\in B\}$. For the rest of this section, we will assume that $N\cdot B\subseteq \partial(A).$

\begin{lem}\label{dim1} Let $(A,*)$ be a finite dimensional commutative associative algebra with the identity $\hat{1}$ such that $\dim A\geq 2$. Let $B$ be a vertex $A$-algebroid that is not a Lie algebra. Assume that 

\vspace{0.2cm}

\noindent (i) there exist $e,f,h\in B$ such that $e_0f=h$, $h_0e=2e$, $h_0f=-2f$, $h_0h=0$, $Span\{e,f,h\}$ is a Lie algebra that is isomorphic to $sl_2$,

\vspace{0.2cm}

\noindent (ii) $A=\mathbb{C}{\hat{1}}\oplus N$ and $N\cdot B\subseteq \partial(A)$. Here, $N$ is a $sl_2$-module.

\vspace{0.2cm}

\noindent If $N$ is one-dimensional, then there exists $a\in A\backslash\{0\}$ such that $a*a=0$. In addition, $A$ is a trivial $sl_2$-module.
\end{lem}
\begin{proof} Since $\dim N=1$, it follows that $N$ and $A$ are trivial $sl_2$-modules. Since $A$ is a 2-dimensional unital commutative associative algebra, we then have that either $A=Span\{\hat{1}, \alpha~|~\alpha*\alpha=\hat{1}\}$ or $A=Span\{\hat{1},\alpha'~|~\alpha'*\alpha'=0\}$. 

\vspace{0.1cm}

\noindent Suppose that $A=Span\{\hat{1}, \alpha~|~\alpha*\alpha=\hat{1}\}$. 
Since
$$\alpha\cdot(\alpha\cdot e)-e=\alpha\cdot (\alpha\cdot e)-(\alpha*\alpha)\cdot e=(e_0\alpha)\cdot \partial (\alpha)+(e_0\alpha)\cdot \partial (\alpha)=0,$$ and $\alpha\cdot (\alpha\cdot e)\in \partial(A)$, we then have that $e\in\partial(A)$. This is a contradiction. Hence, $A=Span\{\hat{1},\alpha'\}$ such that $\alpha'*\alpha'=0$. 
\end{proof}
\begin{lem} Let $(A,*)$ be a finite dimensional commutative associative algebra with the identity $\hat{1}$ such that $\dim A\geq 2$. Let $B$ be a vertex $A$-algebroid that is not a Lie algebra. Assume that 

\vspace{0.2cm}

\noindent (i) there exist $e,f,h\in B$ such that $e_0f=h$, $h_0e=2e$, $h_0f=-2f$, $h_0h=0$, $Span\{e,f,h\}$ is a Lie algebra that is isomorphic to $sl_2$,

\vspace{0.2cm}

\noindent (ii) $A=\mathbb{C}{\hat{1}}\oplus N$ and $N\cdot B\subset \partial(A)$. Here, $N$ is a $sl_2$-module.

\vspace{0.2cm}

\noindent If $e_1f=k\hat{1}$ such that $k\neq 0$ then $\dim N\geq 2$. 
\end{lem}
\begin{proof} If $\dim N=1$, then by Lemma \ref{dim1} there exists $a\in A\backslash\{0\}$ such that $a*a=0$. Note that $\{\hat{1},a\}$ is a basis of $A$. 

We set $a\cdot e=\beta\partial(a)$. Here $\beta\in\mathbb{C}$. Since $\partial(a)_1f=-a_0f=0$, we then have that 
$$ka=a*(e_1f)-e_0f_0a=(a\cdot e)_1 f=\beta\partial (a)_1f=0.$$ This is impossible. Therefore, $\dim N\geq 2$. \end{proof}


\noindent For next Theorem, under suitable conditions, we show that $A$ is a local algebra and obtain a precise description of $A$ as $sl_2$-module.
\begin{thm}\label{sl2insideB} 
Let $(A,*)$ be a finite dimensional commutative associative algebra with the identity $\hat{1}$ such that $\dim A\geq 2$. Let $B$ be a vertex $A$-algebroid that is not a Lie algebra. Assume that there exist $e,f,h\in B$ such that $e_0f=h$, $h_0e=2e$, $h_0f=-2f$, $h_0h=0$, $Span\{e,f,h\}$ is a Lie algebra that is isomorphic to $sl_2$, $A=\mathbb{C}{\hat{1}}\oplus_{i=1}^t N^i$ where $N^i$ are irreducible $sl_2$-modules and $(\oplus_{i=1}^t N^i)\cdot B\subseteq \partial(A)$.

\vspace{0.2cm}

\noindent For each $i\in\{1,...,t\}$, we let $a^i_0$ be the highest weight vector of $N^i$ of weight $m_i$, and set $a^i_j=\frac{1}{j!}(f_0)^ja^i_0$. Note that $\{a^i_0,....,a^i_{m_i}\}$ is a basis of $N^i$ and 
$h_0a^i_j=(m_i-2j)a^i_j$, $f_0a^i_j=(j+1)a^i_{j+1}$, $ e_0a^i_{j}=(m_i-j+1)a^i_{j-1}$.

\vspace{0.2cm} 

\noindent If $e_1f=k\hat{1}$ where $k\neq 0\in\mathbb{C}$, then 

\vspace{0.2cm}

\noindent (i) $Ker(\partial)=\mathbb{C}{\hat{1}}$;

\vspace{0.2cm}

\noindent (ii) each $N^i$ is an irreducible $sl_2$-module that has dimension 2 and $A$ is a local algebra. In addition, for $i,j\in\{1,...,t\}$, $s,r\in\{0,1\}$, the following statements hold: 
\begin{eqnarray*}
&&a^{i}_{s}*a^j_r=0,\\
&&a^i_{0}\cdot e=0,~a^i_{1}\cdot e=\partial(a^i_{0}),\\
&&a^i_{0}\cdot f=\partial(a^i_{1}),~a^i_{1}\cdot f=0,\\
&&a^i_{0}\cdot h=\partial(a^i_{0}),~a^i_{1}\cdot h=-\partial(a^i_{1}),\\
&&a^i_s\cdot\partial(a^j_t)=0,\\
&&\partial(a^i_s)_1e=e_0a^i_{s}=(2-s)a^i_{s-1},\\
&&\partial(a^i_{s})_1f=f_0a^i_{s}=(s+1)a^i_{s+1},~\partial(a^i_{s})_1h=h_0a^i_{s}=(1-2s)a^i_{s},~k=1.
\end{eqnarray*}

\noindent (iii) For $a,a'\in N$, we have $a\cdot e=\partial(e_0a)$, $a\cdot f=\partial(f_0a)$, $a\cdot h=\partial(h_0a)$, $a\cdot \partial(a')=0$. 

\end{thm}

\begin{proof} We will separate the proof of this Theorem into 4 steps. Note we observe that $\partial(\oplus_{i=1}^tN^i)\neq \{0\}$ (because $B$ is not a Lie algebra). 

\vspace{0.2cm}

\noindent {\bf Step 1}: let $i'\in\{1,...,t\}$. We will show that if $\partial(a^{i'}_j)=0$ for some $j\in\{0, 1,...,m_{i'}\}$ then $N^{i'}\subset Ker(\partial)$. 

\vspace{0.2cm}

\noindent {\bf Proof of Step 1} 

\noindent Assume that there exists $j\in\{0, 1,...,m_{i'}\}$ such that $\partial(a^{i'}_j)=0$. Notice that if $m_{i'}=0$ then $N^{i'}\subset Ker(\partial)$. We now suppose that $m_{i'}\geq 1$. Using the fact that for $p\in\{0,...,m_{i'}\}$, we have $f_0a^{i'}_p=(p+1)a^{i'}_{p+1}$, $e_0a^{i'}_p=(m^{i'}-p+1)a^{i'}_{p-1}$ and $\partial(u_0a)=u_0\partial(a)$ for all $u\in B, a\in A$, we can conclude that $N^{i'}\subseteq Ker(\partial)$. This completes the proof of Step 1. 

\vspace{0.2cm}

\noindent {\bf Step 2}: we set $I:=\{j\in\{1,...,t\}~|~\partial(N^j)\neq 0\}$. We will show that $k=m_i=1$, $a_0^i\cdot e=0$, $a_1^i\cdot e=\partial(a_0^i)$ for all $i\in I$. 

\vspace{0.2cm}

\noindent {\bf Proof of Step 2} 

\noindent For $i\in \{1,...,t\}$, we set $a^{i}_j\cdot e=\sum_{r\in I}\sum_{s=0}^{m_r}\beta^{i,j}_{r,s}\partial(a^r_s).$ Here, $\beta^{i,j}_{r,s}\in\mathbb{C}$. Recall that for $a\in A$, $u,v\in B$, we have 
$(a\cdot u)_1v=a*(u_1v)-u_0v_0a.$
Observe that
\begin{eqnarray}
a^i_j*(e_1f)-e_0f_0a^{i}_j&&=ka^{i}_j-e_0(j+1)a^i_{j+1}\nonumber\\
&&=\begin{cases}ka^i_j\text{ if }j=m_i\\
ka^{i}_j-(j+1)(m_i-j)a^i_j\text{ if }j\neq m_i,
\end{cases}\label{compare1}\\
(a^i_j\cdot e)_1f&&=\sum_{r\in I}\sum_{s=0}^{m_r}\beta^{i,j}_{r,s}(\partial(a^r_s))_1f\nonumber\\
&&=-\sum_{r\in I}\sum_{s=0}^{m_r}\beta^{i,j}_{r,s}(a^r_s)_0f\nonumber\\
&&=\sum_{r\in I}\sum_{s=0}^{m_r-1}\beta^{i,j}_{r,s}(s+1)a^r_{s+1}.\label{compare2}
\end{eqnarray}
\noindent If there exists $p\in I$ such that $m_p=0$, we then have that  $ka^p_0=\sum_{r\in I, r\neq p}\sum_{s=1}^{m_r}\beta^{i,j}_{r,s-1}sa^r_{s}.$ This is impossible since $k\neq 0$. Therefore, for $i\in I$, $m_i\geq 1$. 

\vspace{0.2cm}

\noindent By equations (\ref{compare1}), (\ref{compare2}), we have $$ka^{i}_0-(m_i)a^i_0=\sum_{r\in I}\sum_{s=0}^{m_r-1}\beta^{i,0}_{r,s}(s+1)a^r_{s+1}.$$ This implies that 
\begin{equation}\label{kmi}
k=m_i\text{ for all }i\in I.\end{equation} In addition, $\beta^{i,0}_{r,s}=0$ for all $r\in I,~s\in\{0,...,m_r-1\}$, and 
\begin{equation}\label{ai0e}
a^i_0\cdot e=\sum_{r\in I}\beta^{i,0}_{r, m_r}\partial(a^r_{m_r}).
\end{equation}
Since $(a^i_0\cdot e)_1h=a^i_0*(e_1h)-e_0h_0a^i_0=-m_ie_0a^i_0=0$ and 
\begin{eqnarray*}
(a^i_0\cdot e)_1h&&=(\sum_{r\in I}\beta^{i,0}_{r, m_r}\partial(a^r_{m_r}) )_1h\\
&&=\sum_{r\in I}\beta^{i,0}_{r,m_r}h_0(a^r_{m_r})\\
&&=-\sum_{r\in I}\beta^{i,0}_{r,m_r}m_ra^r_{m_r},
\end{eqnarray*}
we then have that $\beta^{i,0}_{r,m_r}=0$ for all $r$. Therfore, \begin{equation}\label{a0e}
a^i_0\cdot e=0\text{ for all }i\in I.\end{equation} 
\vspace{0.1cm} 

\noindent Similarly, by equations (\ref{compare1}), (\ref{compare2}), we have $$ka^i_{m_i}=\sum_{r\in I}\sum_{s=0}^{m_r-1}\beta^{i,m_i}_{r,s}(s+1)a^r_{s+1}.$$ Consequently,
$$
\beta^{i,m_i}_{r,s}=
\begin{cases}0\text{ when }r\neq i,\\
0\text{ when }r=i \text{ and } s< m_i-1,\\
1\text{ when } r=i \text{ and }s=m_i-1~ (\text{by }( \ref{kmi})).
\end{cases}$$
Since 
\begin{eqnarray*}
(a^i_{m_i}\cdot e)_1h&&=(\partial(a^i_{m_i-1})+\beta^{i,m_i}_{i, m_i}\partial (a^i_{m_i}))_1h\nonumber\\
&&=h_0(a^i_{m_i-1})+\beta^{i,m_i}_{i, m_i}h_0(a^i_{m_i})\nonumber\\
&&=(m_i-2(m_i-1))a^i_{m_i-1}+\beta^{i,m_i}_{i,m_i}(m_i-2m_i)a^i_{m_i}\nonumber\\
&&=(-m_i+2)a^i_{m_i-1}+\beta^{i,m_i}_{i,m_i}(-m_i)a^i_{m_i}\text{ and }\\
a^i_{m_i}*(e_1h)-e_0h_0a^i_{m_i}&&=-(m_i-2m_i)e_0a^i_{m_i}\nonumber\\
&&=m_i(m_i-m_i+1)a^i_{m_i-1}\nonumber\\
&&=m_ia_{m_i-1}\text{ for }i\in I,
\end{eqnarray*}
we have $m_ia_{m_i-1}=(-m_i+2)a^i_{m_i-1}+\beta^{i,m_i}_{i,m_i}(-m_i)a^i_{m_i}$. Clearly, $\beta^{i,m_i}_{i,m_i}=0$ and $m_i=-m_i+2$. Hence, we have 
\begin{eqnarray}
m_i&&=1,\text{ and }\\
a^i_{1}\cdot e&&=\partial(a^i_0)\label{a1e}
\end{eqnarray}
for $i\in I$. This complete the proof of Step 2. 

\vspace{0.2cm}

\noindent{\bf Step 3}: we will show that $Ker\partial=\mathbb{C}\hat{1}$ and $I=\{1,...,t\}$.

\vspace{0.2cm}

\noindent {\bf Proof of Step 3} 

\noindent Recall that $B$ acts trivially on $Ker(\partial)$. Also, notice that $\{\partial(a^r_0),\partial(a^r_1)~|~r\in I\}$ is a basis of $\partial(A)$, and 
\begin{eqnarray*}
&&h_0\partial(a^r_0)=\partial(a^r_0), ~h_0\partial(a^r_1)=-\partial(a^r_1).
\end{eqnarray*} 
Let $\alpha\in Ker(\partial)\backslash\{0\}$. For simplicity, we set $\alpha=\lambda\hat{1}+a$ where $\lambda\in\mathbb{C}$ and $a\in\oplus_{i=1}^tN^i$. Note that $h_0a=0$ since $h_0\alpha=0$. Since $a\cdot h\in \partial(A)$ and $h_0(a\cdot h)=(h_0a)\cdot h+a\cdot (h_0h)=0$, we can conclude that $a\cdot h=0$. Since $$(a\cdot h)_1h=a*(h_1h)-h_0h_0a=2a,$$ we can conclude that $a=0$. Hence, $\alpha\in \mathbb{C}\hat{1}$ and $Ker(\partial)=\mathbb{C}\hat{1}$. In addition, $I=\{1,...,t\}$. This completes the proof of Step 3. Also, the statement (i) of Theorem \ref{sl2insideB} holds.

\vspace{0.2cm}

\noindent {\bf Step 4}: we will show that for $i,j\in\{1,...,t\}$, we have 
\begin{eqnarray*}
&&a_0^i\cdot f=\partial(a_1^i),~a_1^i\cdot f=0,~a^i_0\cdot h=\partial(a^i_0),~a^i_1\cdot h=-\partial(a^i_1),\\
&&a^i_s\cdot \partial(a^j_q)=0, ~a^i_s*a^j_q=0 \text{ for all }s,q\in\{0,1\}.
\end{eqnarray*}

\vspace{0.2cm}

\noindent{\bf Proof of Step 4}

\noindent Recall that for $a\in A$, $u,v\in B$, we have $u_0(a\cdot v)=a\cdot u_0v+(u_0a)\cdot v$. So, for $i\in I$, we have  
$h_0(a^i_0\cdot f)=-a^i_0\cdot f$ and $h_0(a^i_1\cdot f)=-3a^i_1\cdot f$. Since 
\begin{eqnarray*}
&&\partial(A)=Span\{\partial(a^r_0),\partial(a^r_1)~|~r\in I\},\\
&&h_0\partial(a^r_0)=\partial(a^r_0), ~h_0\partial(a^r_1)=-\partial(a^r_1),
\end{eqnarray*} and $a^i_0\cdot f, a^i_{1}\cdot f\in\partial (A)$, we can conclude that 
\begin{equation}\label{a1f}
a^i_1\cdot f=0,\end{equation} and $a^i_0\cdot f=\sum_{r}\beta^r\partial(a^r_1)$ where $\beta^r\in\mathbb{C}$.  Since 
\begin{eqnarray*}
&&(a^i_0\cdot f)_1h=(\sum_r\beta^r\partial(a^r_1))_1h=\sum_r\beta^r h_0a^r_1=-\sum_r\beta^ra^r_1,\text{ and }\\
&&a^i_0*(f_1h)-f_0h_0a^i_0=-f_0a^{i}_0=-a^i_1,
\end{eqnarray*}
we can conclude that $\beta^r=\begin{cases}1\text{ if }r=i,\\ 0\text{ otherwise.}\end{cases}$ Moreover, we have
\begin{equation}\label{a0f}
a^i_0\cdot f=\partial(a^i_1).\end{equation}

\vspace{0.2cm}

\noindent Since $h_0(a^i_0\cdot h)=a^i_0\cdot h$ and $h_0(a^i_1\cdot h)=-a^i_1\cdot h$, we then have that 
$a^i_0\cdot h=\sum_{r}\beta^r\partial(a^r_0)$ and $a^i_1\cdot h=\sum_{r}\gamma^r\partial(a^r_1)$. Here, $\beta^r,\gamma^r\in\mathbb{C}$. Since
\begin{eqnarray*}
&&(a^i_0\cdot h)_1f=\sum_r\beta^rf_0a^r_0=\sum_r\beta^ra^r_1,\text{ and }\\
&&a^i_0*(h_1f)-h_0f_0a^i_0=-h_0a^i_1=a^i_1,
\end{eqnarray*} 
it follows that $\beta^r=\begin{cases}1\text{ if }r=i\\ 0\text{ otherwise}\end{cases}$. In addition \begin{equation} a^i_0\cdot h=\partial(a^i_0).\end{equation} Similarly, since 
\begin{eqnarray*}
&&(a^i_1\cdot h)_1e=\sum_r\beta^re_0a^r_1=\sum_r\beta^ra^r_0,\\
&&a^i_1*(h_1e)-h_0e_0a^i_1=-h_0a^i_0=-a^i_0,
\end{eqnarray*} 
we then have \begin{equation}a^i_1\cdot h=-\partial(a^i_1). \end{equation}

\vspace{0.2cm}

\noindent Since 
\begin{eqnarray*}
&&h_0a^i_0\cdot\partial(a^j_0)=2a^i_0\cdot \partial(a^j_0),\\
&&h_0a^i_0\cdot\partial(a^j_1)=0,\\
&&h_0a^i_1\cdot\partial(a^j_0)=0,\\
&&h_0a^i_1\cdot\partial(a^j_1)=-2a^i_0\cdot \partial(a^j_1),\\
\end{eqnarray*} 
we can conclude that \begin{equation}\label{asaq}
a^i_s\cdot \partial(a^j_q)=0\end{equation} for all $s,q\in\{0,1\}$.

\vspace{0.2cm}

\noindent For $i,j\in I$, $h_0(a^i_0*a^j_0)=2a^i_0* a^j_0$ and $h_0(a^i_1*a^{j_1})=-2a^i_1*a^j_1$, we can conclude that \begin{equation}
a^i_0*a^j_0=0\text{ and }a^i_1*a^j_1=0.
\end{equation}
Since $h_0a^i_0*a^j_1=0$, we conclude that there exists $\beta\in \mathbb{C}$ such that $a^i_0*a^j_1=\beta\hat{1}$. Since $0=(a^i_0*a^i_0)*a^j_1=a^i_0*(a^i_0*a^j_1)=a^i_0*\beta\hat{1}=\beta a^i_0$, it implies that $\beta=0$ and 
\begin{equation}\label{aiaj}
a^i_0*a^j_1=0.
\end{equation}
This completes the proof of Step 4. 

\vspace{0.2cm}

\noindent Recall that $N=Span\{a^i_0, a^i_1~|~i\in\{1,....,t\}\}$ and $I=\{1,...t\}$. By equations (\ref{a0e}), (\ref{a1e}), (\ref{a1f})-(\ref{asaq}), we can conclude that 
$a\cdot e=\partial(e_0a)$, $a\cdot f=\partial(f_0a)$, $a\cdot h=\partial(h_0 a)$, $a\cdot \partial(a')=0$ for all $a,a'\in N$. By (\ref{aiaj}), $A$ is a local algebra. This completes the proof of Theorem \ref{sl2insideB}.
\end{proof}


\section{From Lie algebras to Vertex Algebroids} 
 
\noindent In this section, we will first construct vertex algebroids $B_{\mathfrak{g}}$ from unital commutative associative algebras $A$ and Lie algebras $\mathfrak{g}$. Next, in subsection 3.1, we study the algebraic structure of the constructed vertex algebroids $B_{\mathfrak{g}}$ when $\mathfrak{g}$ are finite dimensional simple Lie algebras. Finally, in subsection 3.2, we investigate the algebraic structure of vertex algebroids $B$ when $B$ are simple Leibniz algebras. In particular, under a suitable condition, we show that the vertex algebroids associated with simple Leibniz algebras are equivalent to vertex algebroids that are constructed in subsection 3.1. It is worth to mention that we provide necessary background on representation theory of Leibniz algebras in subsection 5.1.

\vspace{0.2cm}

\noindent Let $A$ be a commutative associative algebra with the identity $\hat{1}$ over $\mathbb{C}$. Let $\mathfrak{g}$ be a Lie algebra that acts on $A$ as a Lie algebra of derivations. Next, we let $U$ be a $\mathfrak{g}$-module such that $\partial:A \rightarrow U$ is a $\mathfrak{g}$-module epimorphism and $\partial(\hat{1})=0$. We set $$B_{\mathfrak{g}}:=\mathfrak{g}+\partial(A).$$ We define a bilinear map $[~,~]$ on $B_{\mathfrak{g}}$ in the following way: for $u,v\in \mathfrak{g}$, $m,n\in \partial(A)$, 
$$[u+m, v+n]=[u,v]+u\cdot n.$$ It is well known that $B_{\mathfrak{g}}$ is a left-Leibniz algebra. 

\vspace{0.2cm} 

\noindent Clearly, $A$ is a module of the left-Leibniz algebra $B_{\mathfrak{g}}$ under the following action: for $u\in\mathfrak{g}, m\in \partial(A)$, $a\in A$,
$$(u+m)\cdot a=u\cdot a.$$

\vspace{0.2cm}

\noindent Let $\langle\cdot,\cdot\rangle$ be a symmetric invariant bilinear form on $\mathfrak{g}$. Following Proposition 21 of [JY1], for $i\in\{0,1\}$, we define a bilinear map $(v,w)\mapsto v_iw$ of degree $-i-1$ on $A\oplus B_{\mathfrak{g}}$ in the following way: for $g,g'\in\mathfrak{g},\lambda,\beta\in B_{\mathfrak{g}}, a,a'\in A$,
\begin{eqnarray*}
&&\lambda_0\beta=[\lambda,\beta],~\lambda_0a=\lambda\cdot a,~a_0\lambda=-\lambda\cdot a, ~a_0a'=0,\\
&&g_1g=\langle g,g'\rangle,~ g_1\partial(a)=g\cdot a, ~\partial(a)_1g=g\cdot a, ~\partial(a)_1\partial(a')=0,\\
&&a_1a'=0,~a_1\lambda=0,~\lambda_1a=0.
\end{eqnarray*}
Then the vector space $A\oplus B_{\mathfrak{g}}$ is a 1-truncated conformal algebra. 

\vspace{0.2cm}

\noindent Now, we assume that 
\begin{enumerate}
\item $A=\mathbb{C}\hat{1}\oplus N$ where $N$ is a $\mathfrak{g}$-module.
\item $\partial$ is a bijection map from $N$ onto $\partial(N)$.
\end{enumerate}
We define a $\mathbb{C}$-bilinear map $A\times B_{\mathfrak{g}}\rightarrow B_{\mathfrak{g}}, (a,v)\mapsto a\cdot v$ in the following way: for $a\in N$, $a'\in A$, $g\in\mathfrak{g}$.
\begin{eqnarray*}
&&1\cdot (g+\partial (a'))=g+\partial(a')\\
&&a\cdot (g+ \partial(a'))=\partial(g_0a)
\end{eqnarray*}

\noindent Now, we will study necessary and sufficient conditions for $B_{\mathfrak{g}}$ to be a vertex $A$-algebroid. 

\vspace{0.2cm}

\noindent Recall that for a 1-truncated conformal algebra $A\oplus B$, $B$ is a vertex $A$-algebroid if and only if it satisfies the following conditions: for $a,a'\in A$, $u,v\in B$, $i=0,1$,
\begin{eqnarray}
\label{VAB1}&&a\cdot(a'\cdot u)-(a*a')\cdot u=(u_0a)\cdot \partial a'+(u_0a')\cdot \partial a,\\
\label{VAB2}&&u_0(a\cdot v)-a\cdot (u_0v)=(u_0a)\cdot v,\\
\label{VAB3}&&u_0(a*a')=a*(u_0a')+(u_0a)*a',\\
\label{VAB4}&&a_0(a'\cdot v)=a'*(a_0v),\\
\label{VAB5}&&(a\cdot u)_1v=a*(u_1v)-u_0v_0a,\\
\label{VAB6}&&\partial(a*a')=a\cdot \partial(a')+a'\cdot \partial(a).
\end{eqnarray}

\noindent In the following Lemma, we will show that the equations (\ref{VAB2}), and (\ref{VAB3}) hold for $B_{\mathfrak{g}}$.
\begin{lem} For $\alpha,\alpha'\in A$, $u,v\in B_{\mathfrak{g}}$, we have 
\begin{eqnarray}
\label{Algebroid1}&&u_0(\alpha *\alpha')=\alpha*(u_0\alpha')+(u_0\alpha)*\alpha',\\
\label{Algebroid2}&&u_0(\alpha\cdot v)-\alpha\cdot (u_0v)=(u_0\alpha)\cdot v.
\end{eqnarray} 
\end{lem} 
\begin{proof} To prove statment (\ref{Algebroid1}), it is enough to show that
$$(g+\partial(a))_0(\alpha*\alpha')=((g+\partial(a))_0\alpha)*\alpha'+\alpha*((g+\partial(a))_0\alpha')$$ 
for all $g\in \mathfrak{g},a, \alpha,\alpha'\in A$. 

\vspace{0.1cm}

\noindent Let $g\in \mathfrak{g},a,\alpha,\alpha'\in A$. By using the fact that the Lie algebra $\mathfrak{g}$ acts as a Lie algebra of derivations on $A$ and $\partial(A)$ acts trivially on $A$, we can easily show that
\begin{eqnarray*}
&&(g+\partial(a))_0(\alpha*\alpha')\\
&&=g_0(\alpha*\alpha')\\
&&=(g_0\alpha)*\alpha'+\alpha*(g_0(\alpha'))\\
&&=((g+\partial(a))_0\alpha)*\alpha'+\alpha*((g+\partial(a))_0\alpha').
\end{eqnarray*}
Hence, the statement (\ref{Algebroid1}) holds.

\vspace{0.2cm}

\noindent Next, we will prove statement (\ref{Algebroid2}). It is enough to show for $\alpha, a,a'\in A$, $g,g'\in \mathfrak{g}$, we have 
\begin{eqnarray*}
&&(g+\partial(a))_0(\alpha\cdot (g'+\partial(a'))-\alpha\cdot((g+\partial(a))_0(g'+\partial(a'))\\
&&=((g+\partial(a))_0\alpha)\cdot (g'+\partial(a')).
\end{eqnarray*}
Clearly, there exist $\beta\in \mathbb{C},\lambda\in N$ such that $\alpha=\beta\hat{1}+\lambda$. By direct calculation, we have
\begin{eqnarray*}
&&(g+\partial(a))_0(\alpha\cdot (g'+\partial(a'))-\alpha\cdot((g+\partial(a))_0(g'+\partial(a'))\\
&&=g_0((\beta\hat{1}+\lambda)\cdot (g'+\partial(a'))-(\beta\hat{1}+\lambda)\cdot g_0(g'+\partial(a'))\\
&&=g_0(\beta g'+\beta\partial(a')+\partial(g'_0\lambda))-(\beta g_0g'+\beta \partial(g_0a')+\partial((g_0g')_0\lambda))\\
&&=\beta g_0g'+\beta \partial(g_0a')+\partial(g_0g'_0\lambda)-(\beta g_0g'+\beta \partial(g_0a')+\partial((g_0g')_0\lambda))\\
&&=\partial(g'_0g_0\lambda)\\
&&=(g_0\lambda)\cdot (g'+\partial(a'))\\
&&=(g_0(\beta\hat{1}+\lambda))\cdot (g'+\partial(a'))\\
&&=((g+\partial(a))_0\alpha)\cdot (g'+\partial(a'))\text{ as desired}.
\end{eqnarray*}
Hence, the equation (\ref{Algebroid2}) holds. 
\end{proof}

\noindent Next, we will study relation (\ref{VAB6}).
\begin{lem} The following statements are equivalent: 

\noindent (i) For every $a,a'\in N$, $a*a'=0$;

\noindent (ii) For every $\alpha,\alpha'\in A$, $\partial(\alpha*\alpha')=\alpha\cdot\partial(\alpha')+\alpha'\cdot\partial(\alpha)$.
\end{lem}
\begin{proof} Let $\alpha,\alpha'\in A$. For convenience, we set $\alpha=\beta\hat{1}+\lambda$, $\alpha'=\beta'\hat{1}+\lambda'$. Here, $\beta,\beta'\in\mathbb{C}$ and $\lambda,\lambda'\in N$. Notice that
\begin{eqnarray*}
\partial(\alpha*\alpha')&&=\partial((\beta\hat{1}+\lambda)*(\beta'\hat{1}+\lambda'))\\
&&=\partial(\beta\beta'\hat{1}+\beta\lambda'+\beta'\lambda+\lambda*\lambda')\\
&&=\beta\partial(\lambda')+\beta'\partial(\lambda)+\partial(\lambda*\lambda')
\end{eqnarray*}
and 
$$\alpha\cdot\partial(\alpha')+\alpha'\cdot\partial(\alpha)=(\beta\hat{1}+\lambda)\cdot\partial(\lambda')+(\beta'\hat{1}+\lambda')\cdot \partial(\lambda)=\beta\partial(\lambda')+\beta'\cdot\partial(\lambda).$$ 

\vspace{0.2cm} 

\noindent Consequently, $\partial(\alpha*\alpha')=\alpha\cdot\partial(\alpha')+\alpha'\cdot\partial(\alpha)$ if and only if $\partial(\lambda*\lambda')=0$. Using the fact that $\partial: N\rightarrow \partial(N)$ is injective, we can conclude that $\partial(\alpha*\alpha')=\alpha\cdot\partial(\alpha')+\alpha'\cdot\partial(\alpha)$ if and only if $\lambda*\lambda'=0$. In addition, the statements (i) and (ii) are equivalent.
\end{proof}

\begin{thm}\label{LiealgebraVAB} Let $A$ be a commutative associative algebra with the identity $\hat{1}$ over $\mathbb{C}$. Let $\mathfrak{g}$ be a Lie algebra that acts on $A$ as a Lie algebra of derivations. Let $U$ be a $\mathfrak{g}$-module such that $\partial:A \rightarrow U$ is a $\mathfrak{g}$-module epimorphism and $\partial(\hat{1})=0$.  We set $B_{\mathfrak{g}}=\mathfrak{g}+\partial(A)$. Then 

\vspace{0.2cm}

\noindent (i) the vector space $(B_{\mathfrak{g}},[~,~])$ is a left Leibniz algebra under the following action: for $u,v\in \mathfrak{g}$, $m,n\in \partial(A)$, 
$$[u+m, v+n]=[u,v]+u\cdot n.$$ In addition, $A$ is a $B_{\mathfrak{g}}$-module under the following action: for $u\in\mathfrak{g}, m\in \partial(A)$, $a\in A$,
$(u+m)\cdot a=u\cdot a.$

\vspace{0.2cm}

\noindent (ii) Let $\langle\cdot,\cdot\rangle$ be a symmetric invariant bilinear form on $\mathfrak{g}$. Then the vector space $A\oplus B_{\mathfrak{g}}$ is a 1-truncated conformal algebra under the following action: for $g,g'\in\mathfrak{g},\lambda,\beta\in B_{\mathfrak{g}}, a,a'\in A$,
\begin{eqnarray*}
&&\lambda_0\beta=[\lambda,\beta],~\lambda_0a=\lambda\cdot a,~a_0\lambda=-\lambda\cdot a, ~a_0a'=0,\\
&&g_1g=\langle g,g'\rangle,~ g_1\partial(a)=g\cdot a, ~\partial(a)_1g=g\cdot a, ~\partial(a)_1\partial(a')=0,\\
&&a_1a'=0,~a_1\lambda=0,~\lambda_1a=0.
\end{eqnarray*}

\noindent (iii) Assume that $A=\mathbb{C}\hat{1}\oplus N$, $\partial$ is a bijection map from $N$ onto $\partial(N)$ and $a*a'=0$ for all $a,a'\in N$. Here $N$ is a $\mathfrak{g}$-module. We define a $\mathbb{C}$-bilinear map $A\times B_{\mathfrak{g}}\rightarrow B_{\mathfrak{g}}, (a,v)\mapsto a\cdot v$ in the following way: for $a\in N$, $a'\in A$, $g\in\mathfrak{g}$.
\begin{eqnarray*}
&&1\cdot (g+\partial (a'))=g+\partial(a'),\\
&&a\cdot (g+ \partial(a'))=\partial(g_0a).
\end{eqnarray*}
Then $B_{\mathfrak{g}}$ is a vertex $A$-algebroid if and only if $g_0g'_0a+g'_0g_0a=(g_1g')*a$ for all $a\in N, g,g'\in\mathfrak{g}$. 
\end{thm}
\begin{proof} The statements (i), (ii) are clear. Now, to show that $B_{\mathfrak{g}}$ is a vertex $A$-algebroid, we need to prove that 
for $\alpha,\alpha'\in A, u,v\in B_{\mathfrak{g}}$, we have
\begin{eqnarray}
\label{Algebroid3}&&\alpha\cdot(\alpha'\cdot u)-(\alpha*\alpha')\cdot u=(u_0\alpha)\cdot \partial(\alpha')+(u_0\alpha')\cdot \partial(\alpha),\\
\label{Algebroid4}&&\alpha_0(\alpha'\cdot v)=\alpha'*(\alpha_0v),\\
\label{Algebroid5}&&(\alpha\cdot u)_1v=\alpha*(u_1v)-u_0v_0\alpha.
\end{eqnarray}

\vspace{0.1cm}

\noindent To show that (\ref{Algebroid3}) holds, it is enough to show that $$\alpha\cdot (\alpha'\cdot (g+\partial(a)))-(\alpha*\alpha')\cdot (g+\partial(a)) =((g+\partial(a))_0\alpha)\cdot \partial(\alpha')+((g+\partial(a))_0\alpha')\cdot \partial(\alpha)$$ for every $\alpha,\alpha' \in A, g\in\mathfrak{g}, a\in N$. Let $\alpha,\alpha'\in A$, $g\in \mathfrak{g}, a\in N$. There exist $\beta,\beta'\in\mathbb{C}$ and $\lambda,\lambda'\in N$ such that $\alpha=\beta\hat{1}+\lambda, \alpha'=\beta'\hat{1}+\lambda'$.  
Notice that 
\begin{eqnarray*}
&&((g+\partial(a))_0\alpha)\cdot \partial(\alpha')+((g+\partial(a))_0\alpha')\cdot \partial(\alpha)\\
&&=(g_0\lambda)\cdot \partial(\lambda')+(g_0\lambda')\cdot \partial(\lambda)\\
&&=0~ (\text{because } g_0\lambda, g_0\lambda'\in N).
\end{eqnarray*}
Since
\begin{eqnarray*}
&&\alpha\cdot(\alpha' \cdot (g+ \partial(a)))-(\alpha*\alpha')\cdot (g+\partial(a))\\
&&=(\beta\hat{1}+\lambda)\cdot((\beta'\hat{1}+\lambda')\cdot (g+\partial(a)))-(\beta\beta'+\beta\lambda'+\beta'\lambda+\lambda*\lambda')\cdot (g+\partial(a))\\
&&=(\beta\hat{1}+\lambda)\cdot(\beta' g+\beta'\partial(a)+\partial(g_0\lambda'))-(\beta\beta' g+\beta\beta'\partial(a)+\beta\partial(g_0\lambda')+\beta'\partial(g_0\lambda))\\
&&=\beta\beta'g+\beta\beta'\partial(a)+\beta\partial(g_0\lambda')+\beta'\partial(g_0\lambda)-(\beta\beta' g+\beta\beta'\partial(a)+\beta\partial(g_0\lambda')+\beta'\partial(g_0\lambda))\\
&&=0,
\end{eqnarray*}
we can conclude immediately that $$\alpha\cdot (\alpha'\cdot (g+\partial(a)))-(\alpha*\alpha')\cdot (g+\partial(a))=((g+\partial(a))_0\alpha)\cdot \partial(\alpha')+((g+\partial(a))_0\alpha')\partial(\alpha).$$ Hence, the statement (\ref{Algebroid3}) holds.

\vspace{0.2cm} 

\noindent Next, we will show that the statement (\ref{Algebroid4}) is true. We only need to show that for $\alpha,\alpha'\in A$, $g\in \mathfrak{g}, a\in N$, $\alpha_0(\alpha'\cdot (g+\partial(a)))=\alpha'*(\alpha_0(g+\partial(a)))$. 

\vspace{0.1cm}

\noindent For simplicity, we set $\alpha=\beta\hat{1}+\lambda$, $\alpha'=\beta'\hat{1}+\lambda'$ where $\beta,\beta'\in\mathbb{C}$, $\lambda,\lambda'\in N$. Using the fact that $N$ is a $\mathfrak{g}$-module, and the assumption that $a*a'=0$ for all $a,a'\in N$, we have 
\begin{eqnarray*}
\alpha_0(\alpha'\cdot (g+\partial(a)))&&=(\beta\hat{1}+\lambda)_0((\beta'+\lambda')\cdot ( g+\partial(a)))\\
&&=\lambda_0(\beta' g+\beta'\partial(a)+\partial(g_0\lambda'))\\
&&=\beta' \lambda_0g\\
&&=(\beta'\hat{1}+\lambda')*(\lambda_0g)\\
&&=\alpha'*(\alpha_0g)\\
&&=\alpha'*(\alpha_0(g+\partial(a))).
\end{eqnarray*}
Hence, the statement (\ref{Algebroid4}) holds.

\vspace{0.2cm}

\noindent Finally, we will study equation (\ref{Algebroid5}). Let $\alpha\in A, u,v\in B_{\mathfrak{g}}$. We will set $\alpha=\beta\hat{1}+\lambda$, $u=g+\partial(a)$, $v=g'+\partial(a')$. Here, $\beta\in\mathbb{C}$, $\lambda,a,a'\in N$. Notice that 
\begin{eqnarray*}
&&\alpha*(u_1v)-u_0v_0\alpha\\
&&=(\beta\hat{1}+\lambda)*((g+\partial(a))_1(g'+\partial(a'))-g_0g'_0\lambda\\
&&=(\beta\hat{1}+\lambda)*(g_1g'+g_0a'+g'_0a)-g_0g'_0\lambda\\
&&=\beta g_1g'+\beta g_0a'+\beta g'_0a+\lambda *(g_1g')-g_0g'_0\lambda\\
\end{eqnarray*}
and 
\begin{eqnarray*}
&&(\alpha\cdot u)_1v\\
&&=((\beta\hat{1}+\lambda)\cdot (g+\partial(a)))_1(g'+\partial(a'))\\
&&=(\beta g+\beta\partial(a)+\partial(g_0\lambda))_1(g'+\partial(a'))\\
&&=\beta g_1g'+\beta g_0a'+\beta g'_0a+ g'_0g_0\lambda.
\end{eqnarray*}
We can conclude immediately that 
$\alpha*(u_1v)-u_0v_0\alpha=(\alpha\cdot u)_1v$ if and only if $g'_0g_0\lambda=\lambda *(g_1g')-g_0g'_0\lambda$. Therefore $B_{\mathfrak{g}}$ is a vertex $A$-algebroid if and only if $g_0g'_0a+g'_0g_0a=(g_1g')*a$ for every $a\in N, g,g'\in\mathfrak{g}$. This completes the proof of this Theorem.
\end{proof}
\begin{ex} Let $\mathfrak{g}=Span\{e,f,h\}$ be a simple Lie algebra that is isomorphic to $sl_2$. Assume that $\mathfrak{g}$ is equipped with a symmetric invariant bilinear form $\langle~,~\rangle$ such that $\langle e,f\rangle\neq 0$. Let $N=Span\{a^0,...,a^m\}$ be an irreducible $sl_2$-module such that $h_0a^i=(m-2i)a^i$, $f_0a^i=(i+1)a^{i+1}$, $e_0a^i=(m-i+1)a^{i-1}$. We set $A=\mathbb{C}\hat{1}\oplus N$. We define a multiplication $*$ on $A$ in the following way: for $\beta,\beta'\in\mathbb{C}$, $a,a'\in N$, $(\beta\hat{1}+a)*(\beta'\hat{1}+a')=(\beta\beta')\hat{1}+\beta a'+\beta' a$. Hence, $A$ is a commutative associative algebra with the identity $\hat{1}$. 

\vspace{0.2cm}

\noindent We set $B=sl_2+\partial(A)$. Assume that $B$ is a vertex $A$-algebroid. Notice that 
$$ma^0=e_0a^1=e_0f_0a^0+f_0e_0a^0=(e_1f)*a^0=\langle e,f\rangle a^0.$$ Hence, we have $\langle e,f\rangle=m$. Since 
$$2(m-1)a^1 =ma^1+(m-2)a^1=f_0h_0a^0+h_0f_0a^0=(h_1f)* a=\langle h,f\rangle a=0,$$ we have $m=1$, $\langle e,f\rangle=1$, and $\dim N=2$. 
\end{ex}

\subsection{Vertex Algebroids associated with simple Lie algebras}

\ \

\vspace{0.2cm}

\noindent In this subsection, we investigate the algebraic structure of vertex algebroids $B_{\mathfrak{g}}$ when $\mathfrak{g}$ is a simple Lie algebra and $N$ is a finite dimensional simple module $L_{\lambda}$ where $\lambda$ is dominant integral.

\vspace{0.2cm}

\noindent  Let $\mathfrak{g}$ be a finite dimensional simple Lie algebra. Let $\langle\cdot,\cdot\rangle:\mathfrak{g}\times \mathfrak{g}\rightarrow\mathbb{C}$ be a nondegenerate symmetric invariant bilinear form, and let $\mathfrak{h}$ be a Cartan subalgebra of $\mathfrak{g}$. We identify $\mathfrak{h}$ with $\mathfrak{h}^*$ by means of the form $\langle\cdot,\cdot\rangle$. Let $\Delta$ be the root system of $\mathfrak{g}$ with respect to $\mathfrak{h}$, viewed as a subset of $\mathfrak{h}$($=\mathfrak{h}^*$), and assume that the form is normalized so that $\langle\alpha,\alpha\rangle=2$ for long roots $\alpha\in\Delta$. We have the root space decomposition 
$$\mathfrak{g}=\mathfrak{h}\oplus \sum_{\alpha\in\Delta}\mathfrak{g}_{\alpha},$$ where $\mathfrak{g}_{\alpha}=\{a\in\mathfrak{g}~|~[h,a]=\alpha(h)a\text{ for }h\in\mathfrak{h}\}$, for $\alpha\in\Delta$. We have $\dim \mathfrak{g}_{\alpha}=1$ and $\dim[\mathfrak{g}_{\alpha},\mathfrak{g}_{-\alpha}]=1$ for $\alpha\in\Delta$. For $\alpha\in \Delta$, we let $h_{\alpha}$ be the unique vector in $\mathfrak{h}$ such that $[\mathfrak{g}_{\alpha},\mathfrak{g}_{-\alpha}]=\mathbb{C}\mathfrak{h}_{\alpha}$ and $\alpha(h_{\alpha})=2$. Then under the identification between $\mathfrak{h}$ and $\mathfrak{h}^*$ we have $h_{\alpha}=\frac{2\alpha}{\langle\alpha,\alpha\rangle},$ and for $a\in\mathfrak{g}_{\alpha}$ and $b\in\mathfrak{g}_{-\alpha}$, we have $[h_{\alpha},a]=2a,~[h_{\alpha},b]=-2b,~[a,b]=\langle a,b\rangle\alpha=\langle a,b\rangle\frac{2}{\langle h_{\alpha},h_{\alpha}\rangle} h_{\alpha};$
we also have $\frac{1}{2}\langle\alpha,\alpha\rangle=\frac{2}{\langle h_{\alpha},h_{\alpha}\rangle}$. We set $n_{\alpha}=\frac{1}{2}\langle h_{\alpha},h_{\alpha}\rangle=\frac{2}{\langle\alpha,\alpha\rangle}$. Then $n_{\alpha}$ is always a positive integer, and in fact $n_{\alpha}=1,2$, or $3$, according to whether $\alpha$ is a long root, a short root if $\mathfrak{g}$ is of type $B_n$, $C_n$ or $F_4$, or a short root if $\mathfrak{g}$ is of type $G_2$. The case in which $\alpha$ is a long root (i.e., $n_{\alpha}=1$) is the case $\langle\alpha,\alpha\rangle=2$, that is, $h_{\alpha}=\alpha$. For $\alpha\in \Delta$, we set
\begin{equation}\mathfrak{g}^{\alpha}=\mathfrak{g}_{\alpha}\oplus \mathbb{C}h_{\alpha}\oplus\mathfrak{g}_{-\alpha}.\end{equation} $\mathfrak{g}^{\alpha}$ is a subalgebra isomorphic to the three-dimensional simple Lie algebra $sl_2$. 

\vspace{0.2cm} 

\noindent Now, we fix a system $\Delta_+\subset \Delta$ of positive roots. We have the triangular decomposition \begin{equation}\mathfrak{g}=\mathfrak{g}_+\oplus\mathfrak{h}\oplus\mathfrak{g}_{-},\end{equation} where $\mathfrak{g}_{\pm}=\sum_{\alpha\in\Delta_+}\mathfrak{g}_{\pm\alpha}.  $ For $\lambda\in\mathfrak{h}^*$, we denote by $V(\lambda)$ the Verma $\mathfrak{g}$-module with highest weight $\lambda$: 
$$V(\lambda)=U(\mathfrak{g})\otimes_{U(\mathfrak{h}\oplus\mathfrak{g}_+)}\mathbb{C}_{\lambda},$$ where $\mathbb{C}_{\lambda}$ is the $\mathfrak{h}\oplus \mathfrak{g}_+$-module on which $\mathfrak{h}$ acts according to the character $\lambda$  and $\mathfrak{g}_+$ acts trivially. We also denote by $L(\lambda)$ the irreducible highest weight $\mathfrak{g}$-module, the quotient of $V(\lambda)$ by its (unique) maximal proper submodule. The $\mathfrak{g}$-module $L(\lambda)$ is finite dimensional if and only if $\lambda$ is dominant integral in the sense that $\lambda(h_{\alpha})=\frac{2\langle\lambda,\alpha\rangle}{\langle\alpha,\alpha\rangle}$ for $\alpha\in \Delta_+$. 

\vspace{0.2cm}

\noindent Let $\theta$ be the highest root, which is also the highest weight of the adjoint $\mathfrak{g}$-module. Notice that $\langle\theta,\theta\rangle=2$. Fix nonzero vectors $e_{\theta}\in\mathfrak{g}_{\theta}$, and $f_{\theta}\in\mathfrak{g}_{-\theta}$ such that $\langle e_{\theta},f_{\theta}\rangle=1$. We have 
$$[h_{\theta},e_{\theta}]=2e_{\theta},~[h_{\theta},f_{\theta}]=-2f_{\theta},~[e_{\theta},f_{\theta}]=h_{\theta}.$$

\noindent By Theorem \ref{LiealgebraVAB}, we have the following proposition. 
\begin{prop}\label{simpleLiealgebraVAB}  Let $\mathfrak{g}$ be a finite dimensional simple Lie algebra equipped with a non degenerate symmetric invariant bilinear form $\langle~,~\rangle$. For a dominant integral $\lambda\in \mathfrak{h}^*$, we let $A_{\lambda}=\mathbb{C}\hat{1}\oplus N_{\lambda}$ be a vector space such that $N_{\lambda}\cong L_{\lambda}$ as $\mathfrak{g}$-modules. We define muliplication $*$ on $A_{\lambda}$ in the following way:
$$\hat{1}*a=a=a*\hat{1},~\alpha*\alpha'=0\text{ for all }a\in A_{\lambda},\alpha,\alpha'\in N_{\lambda}.$$ 
We denote a $\mathfrak{g}$-module isomorphism map from $N_{\lambda}$ to $L_{\lambda}$ by $\partial$. Also, we will extend this map to be a map from $A_{\lambda}$ to $L_{\lambda}$ by assuming map $\partial(\hat{1})=0$. 

\noindent Now, we set $$B_{\lambda}=\mathfrak{g}+\partial(A_{\lambda}).$$ Then 

\vspace{0.2cm}

\noindent (i) the vector space $(B_{\lambda},[~,~])$ is a left Leibniz algebra under the following action: for $u,v\in \mathfrak{g}$, $m,n\in \partial(A_{\lambda})$, 
$$[u+m, v+n]:=[u,v]+u\cdot n.$$ Also, $A_{\lambda}$ is a $B_{\lambda}$-module under the following action: for $u\in\mathfrak{g}, m\in \partial(A_{\lambda})$, $a\in A_{\lambda}$,
$(u+m)\cdot a:=u\cdot a.$

\vspace{0.2cm}

\noindent (ii) Let $\langle\cdot,\cdot\rangle$ be a symmetric invariant bilinear form on $\mathfrak{g}$. Then the vector space $A_{\lambda}\oplus B_{\lambda}$ is a 1-truncated conformal algebra under the following action: for $g,g'\in\mathfrak{g},u,v\in B_{\lambda}, a,a'\in A_{\lambda}$,
\begin{eqnarray*}
&&u_0v:=[u,v],~ u_0a:= u\cdot a,~a_0u:=-u\cdot a, ~a_0a':=0,\\
&&g_1g:=\langle g,g'\rangle,~ g_1\partial(a):=g\cdot a, ~\partial(a)_1g:=g\cdot a, ~\partial(a)_1\partial(a'):=0,\\
&&a_1a':=0,~a_1u:=0,~u_1a:=0.
\end{eqnarray*}

\vspace{0.2cm}

\noindent (iii) We define a $\mathbb{C}$-bilinear map $A_{\lambda}\times B_{\lambda}\rightarrow B_{\lambda}, (a,v)\mapsto a\cdot v$ in the following way: for $a\in N_{\lambda}$, $a'\in A_{\lambda}$, $g\in\mathfrak{g}$,
\begin{eqnarray*}
&&1\cdot (g+\partial (a'))=g+\partial(a'),\\
&&a\cdot (g+ \partial(a'))=\partial(g_0a).
\end{eqnarray*}
Assume that $g_0g'_0a+g'_0g_0a=(g_1g')*a$ for every $a\in N, g,g'\in\mathfrak{g}$. Then $B_{\lambda}$ is a vertex $A_{\lambda}$-algebroid.

\end{prop}


\begin{lem}\label{dominantintegral} For a dominant integral $\lambda\in \mathfrak{h}^*$, we let $A_{\lambda}=\mathbb{C}\hat{1}\oplus N_{\lambda}$ be a vector space such that $N_{\lambda}\cong L_{\lambda}$ as $\mathfrak{g}$-modules. We set $B_{\lambda}=\mathfrak{g}+\partial(A_{\lambda}).$ Assume that $A_{\lambda}$ is a unital commutative associative algebra and $B_{\lambda}$ is a vertex $A_{\lambda}$-algebroid that satisfy Proposition \ref{simpleLiealgebraVAB}. Then
\begin{enumerate}
\item $\mathfrak{g}$ is of $ADE$-type. However, $\mathfrak{g}\neq E_8$.
\item For $\alpha\in \Delta_+$, we have $\lambda(h_{\alpha})=1$. In particular, we have $\lambda(h_{\theta})=1$.
\item We let $v_{\lambda}\in N_{\lambda}$ such that $h\cdot v_{\lambda}=\lambda(h)v_{\lambda}$, $g\cdot v_{\lambda}=0$ for all $h\in \mathfrak{h},g\in\mathfrak{g}_+$. For $\alpha\in \Delta_+$, $a\in N_{\lambda}$, we have 
\begin{eqnarray*}
&&(f_{\alpha})_0(f_{\alpha})_0a=0\text{ and }(e_{\alpha})_0(e_{\alpha})_0a=0. 
\end{eqnarray*}
\end{enumerate}
\end{lem}
\begin{proof} Let $\alpha\in\Delta_+$. Since $\mathfrak{g}_{\alpha}$ is isomorphic to $sl_2$, and $\mathfrak{g}_{\alpha}\subset B_{\lambda}$, by Theorem \ref{sl2insideB}, we can conclude that $(e_{\alpha})_1f_{\alpha}=1$, and $\langle\alpha,\alpha\rangle=(h_{\alpha})_1h_{\alpha}=2$. Therefore, $\mathfrak{g}$ is of $ADE$-type. 

\vspace{0.1cm}

\noindent Recall that 
\begin{equation}\label{algebroidrel}
g_0g'_0a+g'_0g_0a=(g_1g')*a\text{ for every }a\in N_{\lambda},~g,g'\in \mathfrak{g}.\end{equation} Let $\alpha\in\Delta_+$. Since
\begin{eqnarray}
&&(e_{\alpha})_0(f_{\alpha})_0v_{\lambda}+(f_{\alpha})_0(e_{\alpha})_0v_{\lambda}=(e_{\alpha})_0(f_{\alpha})_0v_{\lambda}=(h_{\alpha})_0v_{\lambda}=\lambda(h_{\alpha})v_{\lambda},\text{ and }\\
&&((e_{\alpha})_1f_{\alpha})v_{\lambda}=v_{\lambda},
\end{eqnarray}
we then have $\lambda(h_{\alpha})=1$. In particular, we have $\lambda(h_{\theta})=1$. We can then conclude further that $\mathfrak{g}\neq E_8$. This completes the proofs of statements (1) and (2).

\vspace{0.1cm}

\noindent Since $2(f_{\alpha})_0(f_{\alpha})_0a=((f_{\alpha})_1f_{\alpha})*a=0$ for $a\in N_{\lambda},\alpha\in \Delta_+$, we can conclude that $$(f_{\alpha})_0(f_{\alpha})_0a=0\text{ for all }\alpha\in \Delta_+,a\in N_{\lambda}.$$ Similarly, $$(e_{\alpha})_0(e_{\alpha})_0a=0\text{ for all }\alpha\in \Delta_+,a\in N_{\lambda}.$$ This completes the proof of statement (3). 
\end{proof}

\subsection{Vertex Algebroids associated with Simple Leibniz algebras}

\ \ 

\vspace{0.2cm}

\noindent In this subsection, we investigate the algebraic structure of vertex $A$-algbroids $B$ when $B$ are simple Leibniz algebras. Also, we explore connections between these vertex algebroids and the vertex algebroids that were studied in subsection 3.1. Under a suitable condition, we will show they are equivalent. Note that we provide background information about simple Leibniz algebras in subsection 5.1.

\vspace{0.2cm}

\noindent We let $(A,*)$ be a finite dimensional commutative associative algebra with the identity $\hat{1}$, and let $B$ be a finite dimensional vertex $A$-algebroid such that $B$ is a simple Leibniz algebra that is not a Lie algebra. Then there exists a simple Lie algebra $\mathfrak{g}$ such that $Leib(B)$ is an irreducible module over $\mathfrak{g}$ and $B=\mathfrak{g}+Leib(B)$. Here, $Leib(B)=Span\{u_0u~|~u\in B\}=Span\{u_0v+v_0u~|~u,v\in B\}$. Since $\partial(A)$ is a $\mathfrak{g}$-submodule of $Leib(B)$. We can conclude that $Leib(B)=\partial(A)$. In addition $B=\mathfrak{g}+\partial(A)$. 

\begin{lem} Let $(A,*)$ be a finite dimensional commutative associative algebra with the identity $\hat{1}$, and let $B$ be a finite dimensional vertex $A$-algebroid such that $B$ is a simple Leibniz algebra that is not a Lie algebra. We may and shall set $B=\mathfrak{g}+\partial(A)$. Here $\mathfrak{g}$ is a simple Lie algebra.  Also, we let $A=\mathbb{C}\hat{1}\oplus N$ where $N$ is a $\mathfrak{g}$-submodule of $A$. If $Ker(\partial)=\mathbb{C}\hat{1}$, then $N$ is a simple $\mathfrak{g}$-module.
\end{lem}
\begin{proof} Since $A$ is a module for the simple Lie algebra $\mathfrak{g}$ and $\mathbb{C}\hat{1}$ is a $\mathfrak{g}$-submodule of $A$, we may and shall write $A=\mathbb{C}\hat{1}\oplus N$ where $N$ is a $\mathfrak{g}$-submodule of $A$. Since $\partial(\hat{1})=0$, $Leib(B)=\partial(A)$, $A=\mathbb{C}\hat{1}\oplus N$, we then have $Lieb(B)=\partial(N)$. Clearly, $\partial|_N:N\rightarrow \partial(A)=Leib(B)$ is a $\mathfrak{g}$-module epimorphism. Since $Ker(\partial|_N)=\{u\in N~|~\partial(u)=0\}$ is a $\mathfrak{g}$-submodule of $N$ and $Ker(\partial |_N)\subset Ker(\partial)=\mathbb{C}\hat{1}$, we can conclude that $Ker(\partial |_N)=\{0\}$, and $\partial |_N$ is a $\mathfrak{g}$-module isomorphism. Morever, $N$ is a simple $\mathfrak{g}$-module. 
\end{proof}

\vspace{0.2cm}

\begin{rmk} Recall that $N\cdot B=Span\{u\cdot b~|~u\in N, ~b\in B\}$. Because $v_0(u\cdot b)=(v_0u)\cdot b+u\cdot v_0b$ for all $v\in B$, $u\in N$, $b\in B$, $N\cdot B$ is a left ideal of the Leibniz algebra $B$. Hence $N\cdot B$ is either $\{0\}$ or $\partial(A)$ or $B$. 
\end{rmk}

\noindent Now, we assume that $Ker(\partial)=\mathbb{C}\hat{1}$, and there exist $u,v\in\mathfrak{g}$ such that $u_1v\neq 0$. Since $\mathfrak{g}$ is a Lie algebra, this implies that for $g,g'\in\mathfrak{g}$, $$\partial(g_1g')=0\text{ and }g_1g'\in\mathbb{C}\hat{1}.$$ Let $(~,~):\mathfrak{g}\times\mathfrak{g}\rightarrow\mathbb{C}$ be a bilinear form such that for $g,g'\in \mathfrak{g}$, $$(g,g')\hat{1}=g_1g'.$$ Since $g_1g'=g'_1g$, we can conclude that $(~,~)$ is symmetric. Recall that for $b,b',b''\in B$, we have 
$b''_1(b_0b')=b_0(b''_1b')-(b_0b'')_1b'$. Let $x,y,z\in\mathfrak{g}$. Since $x_1z\in\mathbb{C}\hat{1}$, this implies that 
$$x_1(y_0z)=y_0(x_1z)-(y_0x)_1z=-(y_0x)_1z=(x_0y)_1z.$$ Consequently, $(x,y_0z)=(x_0y,z)$ and $(~,~)$ is invariant. 

\vspace{0.2cm}

\noindent We set $T=\{g\in\mathfrak{g}~|~(g,t)=0\text{ for all }t\in\mathfrak{g}\}$. Observe that for $g\in T$, $b\in\mathfrak{g}$, we have $(t,b_0g)\hat{1}=t_1(b_0g)=b_0(t_1g)-(b_0t)_1g=0$ for all $t\in \mathfrak{g}$. Hence, $b_0g\in T$ and $T$ is an ideal of $\mathfrak{g}$. Since $\mathfrak{g}$ is simple, it follows that $T$ is either $\{0\}$ or $\mathfrak{g}$. Since $(u,v)\hat{1}=u_1v\neq 0$, this implies that $u$ and $v$ are not members of $T$ and $T=\{0\}$. Consequently, $(~,~)$ is non-degenerate symmetric invariant bilinear form on $\mathfrak{g}$ and $(~,~)$ is a scalar multiple of the Killing form of $\mathfrak{g}$. 

\vspace{0.2cm}

\noindent We summarize the above discussion in the following lemma.
\begin{lem}\label{nondegensyminv} Let $(A,*)$ be a finite dimensional commutative associative algebra with the identity $\hat{1}$, and let $B$ be a finite dimensional vertex $A$-algebroid such that $B$ is a simple Leibniz algebra that is not a Lie algebra. Then 

\vspace{0.2cm}

\noindent (i) there exists a simple Lie algebra $\mathfrak{g}$ such that $Leib(B)$ is an irreducible module over $\mathfrak{g}$ and $B=\mathfrak{g}+Leib(B)$. In addition, $Leib(B)=\partial(A)$. 

\vspace{0.2cm}

\noindent (ii) Let $(~,~):\mathfrak{g}\times\mathfrak{g}\rightarrow\mathbb{C}$ be a bilinear form such that $(b,b')\hat{1}=b_1b'$ for $b,b'\in\mathfrak{g}$. If $Ker(\partial)=\mathbb{C}\hat{1}$, and there exist $u,v\in\mathfrak{g}$ such that $u_1v\neq 0$, then $(~,~)$ is a non-degenerate invariant bilinear form, and there exists $\lambda\in \mathbb{C}^*$ such that $(~,~)=\lambda\langle~,~\rangle$. Here $\langle~,~\rangle$ is a Killing form of $\mathfrak{g}$. 
\end{lem}

\vspace{0.2cm}

\noindent Next, we will study the relationship between the vertex algebroids constructed in subsection 4.1 and the vertex $A$-algebroid $B$ that we currently discuss. 

\vspace{0.2cm}

\noindent Following subsection 4.1, we let $\mathfrak{h}$ be a Cartan subalgebra of $\mathfrak{g}$. We identify $\mathfrak{h}$ with $\mathfrak{h}^*$ by means of the form $\langle~,~\rangle$. Let $\Delta$ be the root system of $\mathfrak{g}$ with respect to $\mathfrak{h}$, viewed as a subset of $\mathfrak{h}$ ($=\mathfrak{h}^*$). Note that $\langle\alpha,\alpha\rangle=2$ for long roots $\alpha\in\Delta$. We have the root space decomposition $\mathfrak{g}=\mathfrak{h}\oplus \sum_{\alpha\in\Delta}\mathfrak{g}_{\alpha}$. For $\alpha\in \Delta$, we let $h_{\alpha}$ be the unique vector in $\mathfrak{h}$ such that $[\mathfrak{g}_{\alpha},\mathfrak{g}_{-\alpha}]=\mathbb{C}\mathfrak{h}_{\alpha}$ and $\alpha(h_{\alpha})=2$. For $\alpha\in \Delta$, we set
$\mathfrak{g}^{\alpha}=\mathfrak{g}_{\alpha}\oplus \mathbb{C}h_{\alpha}\oplus\mathfrak{g}_{-\alpha}$ ($\cong sl_2$ as Lie algebras). Also, we let $\theta$ be the highest root, which is also the highest weight of the adjoint $\mathfrak{g}$-module. Therefore, $\langle\theta,\theta\rangle=2$. Fix nonzero vectors $e_{\theta}\in\mathfrak{g}_{\theta}$, and $f_{\theta}\in\mathfrak{g}_{-\theta}$ such that $\langle e_{\theta},f_{\theta}\rangle=1$. We have $[h_{\theta},e_{\theta}]=2e_{\theta},~[h_{\theta},f_{\theta}]=-2f_{\theta},~[e_{\theta},f_{\theta}]=h_{\theta}.$

\begin{thm}\label{vertexalgebroidsimpleLeibnizalgebra} Let $(A,*)$ be a finite dimensional commutative associative algebra with the identity $\hat{1}$, and let $B$ be a finite dimensional vertex $A$-algebroid such that $B$ is a simple Leibniz algebra that is not a Lie algebra. Then there exists a simple Lie algebra $\mathfrak{g}$ such that $B=\mathfrak{g}+\partial(A)$ and $\partial(A)$ is a simple $\mathfrak{g}$-module. In addition, there exists a $\mathfrak{g}$-module $N$ such that $A=\mathbb{C}\hat{1}\oplus N$.

\vspace{0.2cm}

\noindent Let $(~,~):\mathfrak{g}\times\mathfrak{g}\rightarrow\mathbb{C}$ be a bilinear form such that $(b,b')\hat{1}=b_1b'$ for $b,b'\in\mathfrak{g}$. Assume that there exist $u,v\in\mathfrak{g}$ such that $u_1v\neq 0$, $Ker(\partial)=\mathbb{C}\hat{1}$ and $N\cdot B\subset\partial(A)$. Then
\vspace{0.2cm}

\noindent (i) $(~,~)$ is the Killing form of $\mathfrak{g}$. 

\vspace{0.2cm}

\noindent (ii) $A$ is a local algebra. Also, for $a,a'\in N$, $a*a'=0$.

\vspace{0.2cm}

\noindent (iii) For $a\in N$, $u,v\in \mathfrak{g}, a'\in A$, 
\begin{eqnarray*}
&&a\cdot (u+\partial(a'))=\partial(u_0a),\\
&&u_0v_0a+v_0u_0a=(u_1v)*a.
\end{eqnarray*} 
Moreover, there exists a dominant integral $\lambda$ such that $B\cong B_{\lambda}$ as vertex algebroids. 
\end{thm}
\begin{proof} By Lemma \ref{nondegensyminv}, there exists $k\in \mathbb{C}^*$ such that $(~,~)=k\langle~,~\rangle$. Hence $(e_{\theta},f_{\theta})\neq 0$, and $(e_{\theta})_1f_{\theta}\in \mathbb{C}^*\hat{1}$. By Theorem \ref{sl2insideB}, we can conclude that $(e_{\theta})_1f_{\theta}=\hat{1}$, and $k=1$. Therefore, $(~,~)=\langle~,~\rangle$ is the Killing form on $\mathfrak{g}$. In addition, $A$ is a local algebra and $a*a'=0$, $a\cdot\partial(a')=0$ for all $a,a'\in N$. 

\vspace{0.2cm}

\noindent For $\alpha\in \Delta_+$, we fix nonzero vectors $e_{\alpha}\in \mathfrak{g}_{\alpha}$, $f_{\alpha}\in \mathfrak{g}_{-\alpha}$ such that $(h_{\alpha})_0e_{\alpha}=2e_{\alpha}$, $(h_{\alpha})_0f_{\alpha}=-2f_{\alpha}$, $(e_{\theta})_0f_{\alpha}=h_{\alpha}$. Since $\mathfrak{g}_{\alpha}\subseteq B$, $(e_{\alpha},f_{\alpha})=\langle e_{\alpha},f_{\alpha}\rangle\neq 0$, by Theorem \ref{sl2insideB}, we can conclude that $\langle e_{\alpha},f_{\alpha}\rangle= (e_{\alpha},f_{\alpha})=1$ and $a\cdot u=\partial (u_0a)$ for all $a\in N, ~u\in\mathfrak{g}_{\alpha}$. In addition, we have $({\alpha},{\alpha})=2$ for all $\alpha\in\Delta_+$. Consequently, $\mathfrak{g}$ is of $ADE$-type and $a\cdot u=\partial (u_0a)$ for all $a\in N$, $u\in\mathfrak{g}$.

\vspace{0.2cm}

\noindent Let $a\in N$, $u,v\in \mathfrak{g}$. Because $a*(u_1v)-u_0v_0a=(a\cdot u)_1v$ and $a\cdot u=\partial(u_0a)$, we have 
$$a*(u_1v)-u_0v_0a=(\partial(u_0a))_1v=-(u_0a)_0v=v_0u_0a.$$
Because $Ker(\partial)=\mathbb{C}\hat{1}$, we then have that $N$ is a simple $\mathfrak{g}$-module. Hence, there exists an integral domain $\lambda\in\mathfrak{h}^*$ such that $N\cong L_{\lambda}$ as $\mathfrak{g}$-module and $A\cong A_{\lambda}$ as unital commutative associative algebras.

\vspace{0.2cm}

\noindent Now, we recall the definition of vertex algebroid homomorphism (cf. \cite{LiY2}). Let $\mathcal{A}$, $\mathcal{A}'$ be unital commutative associative algebras, and let $\mathcal{B}$ be a vertex $\mathcal{A}$-algebroid, $\mathcal{B}'$ be a vertex $\mathcal{A}'$-algebroid. A {\em{vertex algebroid homomorphism from $\mathcal{B}$ to $\mathcal{B}'$}} is a linear map $\tau:\mathcal{A}\oplus\mathcal{B}\rightarrow\mathcal{A}'\oplus\mathcal{B}'$ such that $\tau(\mathcal{A})\subseteq \mathcal{A}'$, $\tau(\mathcal{B})\subseteq \mathcal{B}'$ and

\noindent (i) $\tau |_{\mathcal{A}}$ is an associative algebra homomorphism;

\noindent (ii) $\tau |_{\mathcal{B}}$ is a Leibniz algebra homomorphism;

\noindent (iii) $\tau(a\cdot b)=\tau(a)\cdot\tau(b)$ for $a\in\mathcal{A}$, $b\in\mathcal{B}$;

\noindent (iv) $\tau(u)_1\tau(v)=\tau(u_1v)$ for $u,v\in\mathcal{B}$;

\noindent (v) $\tau\circ\partial=\partial\circ \tau$;

\noindent (vi) $\tau(b_0a)=\tau(b)_0\tau(a)$ for $a\in\mathcal{A}$, $b\in\mathcal{B}$. 

\vspace{0.1cm}

\noindent Next, we will show that $B$ is isomorphic to $B_{\lambda}$ as vertex algebroids. Recall that 

\noindent (i) $A=\mathbb{C}\hat{1}\oplus N$, $B=\mathfrak{g}+\partial(A)$ and 

\noindent (ii) $A_{\lambda}=\mathbb{C}{\bf 1}\oplus N_{\lambda}$, $B_{\lambda}=\mathfrak{g}+\partial(A_{\lambda})$. 

\noindent Here, $N\cong L_{\lambda}\cong N_{\lambda}$ where $L_{\lambda}$ is an irreducible $\mathfrak{g}$-module. $\hat{1}$ is the identity of the associative commutative algebra $A$ and ${\bf 1}$ is the identity of the associative commutative algebra $A_{\lambda}$. 
 
\noindent Let $\phi:N\rightarrow N_{\lambda}$ be a $\mathfrak{g}$-module isomorphism from $N$ to $N_{\lambda}$. Note that $\partial(A)=\partial(N)$ and $\partial(A_{\lambda})=\partial(N_{\lambda})$ since $\partial(\hat{1})=0$ and $\partial({\bf 1})=0$. Let $\mu:A\oplus B\rightarrow A_{\lambda}\oplus B_{\lambda}$ be a linear map defined by $\mu(\alpha\hat{1}+a)=\alpha{\bf 1}+\phi(a)$ and $\mu(g+\partial(a))=g+\partial(\phi(a))$ where $\alpha\in \mathbb{C}, a\in N$. Clearly, 
$$\mu\partial(a)=\partial(\phi(a))= \partial(\mu(a))\text{ for all }a\in N, $$ and $\mu |_A$ and $\mu |B$ are linear isomorphisms. For $a^1,a^2\in N$, $\alpha^1,\alpha^2\in\mathbb{C}$, we have 
\begin{eqnarray*}
&&\mu |_{A}((\alpha^1\hat{1}+a^1)*(\alpha^2\hat{1}+a^2))\\
&&=\mu |_{A} ((\alpha^1\alpha^2)\hat{1}+\alpha^1 a^2+\alpha^2 a^1)\\
&&=(\alpha^1\alpha^2){\bf 1}+\phi(\alpha^1 a^2)+\phi(\alpha^2 a^1)\\
&&=(\alpha^1\alpha^2){\bf 1}+\alpha^1\phi(a^2)+\alpha^2\phi(a^1)\\
&&=(\alpha^1{\bf 1}+\phi(a^1))*(\alpha^2{\bf 1}+\phi(a^2))\\
&&=\mu |_A(\alpha^1\hat{1}+a^1)*\mu |_{A}(\alpha^2\hat{1}+a^2).
\end{eqnarray*}
Hence, $\mu |_A$ is an algebra homomorphism. 
Next, we will show that $\mu|_{B}$ is a Leibniz algebra homomorphism. Notice that for $g^1,g^2\in\mathfrak{g},a^1,a^2\in N$, we have 
\begin{eqnarray*}
&&\mu((g^1+\partial(a^1))_0(g^2+\partial(a^2))\\
&&=\mu(g^1_0g^2+\partial(g^1_0a^2))\\
&&=g^1_0g^2+\partial(\phi(g^1_0a^2))\\
&&=(g^1+\partial(\phi(a^1)))_0(g^2+\partial(\phi(a^2))\\
&&=\mu(g^1+\partial(a^1))_0\mu(g^2+\partial(a^2)).
\end{eqnarray*}
Hence $\mu|_B$ is a Leibniz algebra homomorphism. Now, we will show that $\mu(a\cdot b)=\mu(a)\cdot \mu(b)$, $\mu(b_0a)=\mu(b)_0\mu(a)$ for  $a\in A$, $b\in B$. Since 
\begin{eqnarray*}
&&\mu((\alpha\hat{1}+a^1)\cdot (g+\partial(a^2))\\
&&=\mu(\alpha g+\alpha\partial(a^2)+\partial(g_0a^1))\\
&&=\alpha g+\alpha\partial(\phi(a^2))+\partial(\phi(g_0a^1))\\
&&=\alpha g+\alpha\partial(\phi(a^2))+\partial(g_0\phi(a^1))\\
&&=(\alpha{\bf 1}+\phi(a^1))\cdot (g+\partial(\phi(a^2))\\
&&=\mu(\alpha\hat{1}+a^1)\cdot\mu(g+\partial(a^2))
\end{eqnarray*}
and
\begin{eqnarray*}
&&\mu((g+\partial(a^2))_0(\alpha\hat{1}+a^1))\\
&&=\mu(g_0 a^1)\\
&&=\phi(g_0a^1)\\
&&=g_0\phi(a^1)\\
&&=(g+\partial(\phi(a^2))_0(\alpha{\bf 1}+\phi(a^1))\\
&&=\mu(g+\partial(a^2))_0\mu(\alpha\hat{1}+a^1)
\end{eqnarray*}
for all $a^1,a^2\in N$, $g\in \mathfrak{g}$, $\alpha\in\mathbb{C}$, we can conclude that $\mu(a\cdot b)=\mu(a)\cdot \mu(b)$, $\mu(b_0a)=\mu(b)_0\mu(a)$ for all $a\in A$, $b\in B$. Finally, we will show that $\mu(u)_1\mu(v)=\mu(u_1v)$ for $u,v\in B$. Observe that 
\begin{eqnarray*}
&&\mu((g^1+\partial(a^1))_1(g^2+\partial(a^2)))\\
&&=\mu(\langle g^1,g^2\rangle\hat{1}+g^1_1\partial(a^2)+\partial(a^1)_1g^2+\partial(a^1)_1\partial(a^2) )\\
&&=\mu(\langle g^1,g^2\rangle\hat{1}-a^2_0g^1-a^1_0g^2-a^1_0\partial(a^2))\\
&&=\langle g^1,g^2\rangle{\bf 1}+\phi(g^1_0a^2)+\phi(g^2_0a^1)\\
&&=\langle g^1,g^2\rangle{\bf 1}+g^1_0\phi(a^2)+g^2_0\phi(a^1)\\
&&=\langle g^1,g^2\rangle{\bf 1}+g^1\cdot\phi(a^2)+g^2\cdot\phi(a^1)\\
&&=(g^1+\partial(\phi(a^1)))_1(g^2+\partial(\phi(a^2)))\\
&&=\mu(g^1+\partial(a^1))_1\mu(g^2+\partial(a^2)).
\end{eqnarray*}
Hence, $\mu(u)_1\mu(v)=\mu(u_1v)$ for $u,v\in B$. Moreover, we can conclude that $B\cong B_{\lambda}$ as vertex algebroids. 
\end{proof}

\section{Vertex Algebra associated with Simple Leibniz algebra}

\noindent In this section, we construct indecomposable non-simple $C_2$-cofinite $\mathbb{N}$-graded vertex algebras from vertex algebroids that were constructed in subsection 3.1. Also, we study modules and conformal vectors of these $\mathbb{N}$-graded vertex algebras. Precisely, in subsection 4.1, we review construction of vertex algebras associated with vertex algebroids and their irreducible $\mathbb{N}$-graded modules. In subsection 4.2, we construct indecomposable non-simple $C_2$-cofinite $\mathbb{N}$-graded vertex algebras associated with simple Leibniz algebras and classify their irreducible $\mathbb{N}$-graded modules. In subsection 4.3, we study conformal vectors of these vertex algebras.

\vspace{0.2cm}

\noindent Background on general theory of vertex algebra is provided in subsection 5.2.

\ \ \ 

\subsection{$\mathbb{N}$-graded Vertex Algebras associated with Vertex Algebroids}

\ \ \ 

\vspace{0.1cm}

\noindent In this subsection, we recall a construction of vertex algebras associated with vertex algebroids in \cite{LiY}. 

\vspace{0.2cm}

\noindent Let $A$ be a commutative associative algebra with identity $\mathfrak{e}$ and let $B$ be a vertex $A$-algebroid. We set $L(A\oplus B)=(A\oplus B)\otimes \mathbb{C}[t,t^{-1}].$ 
Subspaces $L(A)$ and $L(B)$ of $L(A\oplus B)$ are defined in the obvious way. 

\vspace{0.2cm}

\noindent We set $\hat{\partial}=\partial\otimes 1+1\otimes \frac{d}{dt}:L(A)\rightarrow L(A\otimes B).$ We define 
$deg(a\otimes t^n)=-n-1$, $deg(b\otimes t^n)=-n$ for $a\in A, ~b\in B,~n\in\mathbb{Z}$.
Then $L(A\oplus B)$ becomes a $\mathbb{Z}$-graded vector space: 
$$L(A\oplus B)=\oplus_{n\in\mathbb{Z}}L(A\oplus B)_{(n)}$$ where 
$L(A\oplus B)_{(n)}=A\otimes \mathbb{C}t^{-n-1}+B\otimes \mathbb{C}t^{-n}$. Clearly, the subspaces $L(A)$ and $L(B)$ are $\mathbb{Z}$-graded vector spaces as well. In addition, for $n\in \mathbb{N}$, $L(A)_{(n)}=A\otimes  \mathbb{C}t^{-n-1}.$ 
The linear map $\hat{\partial}:L(A)\rightarrow L(A\oplus B)$ is of degree 1. We define a bilinear product $[\cdot,\cdot]$ on $L(A\oplus B)$ as follow:
\begin{eqnarray}
&&[a\otimes t^m,a'\otimes t^n]=0,\label{aa'}\\
&&[a\otimes t^m, b\otimes t^n]=a_0b\otimes t^{m+n},\\
&&[b\otimes t^n,a\otimes t^m]=b_0a\otimes t^{m+n},\\
&&[b\otimes t^m,b'\otimes t^n]=b_0b'\otimes t^{m+n}+m(b_1b')\otimes t^{m+n-1}\label{bb'}
\end{eqnarray}
for $a,a'\in A$, $b,b'\in B$, $m,n\in\mathbb{Z}$. For convenience, we set 
$$\mathcal{L}:=L(A\oplus B)/\hat{\partial}L(A).$$ It was shown in \cite{LiY} that $\mathcal{L}=\oplus_{n\in\mathbb{Z}}\mathcal{L}_{(n)}$ is a $\mathbb{Z}$-graded Lie algebra. Here, 
$$\mathcal{L}_{(n)}=L(A\oplus B)_{(n)}/\hat{\partial}(L(A)_{(n-1)})=(A\otimes \mathbb{C}t^{-n-1}+B\otimes \mathbb{C}t^{-n})/\hat{\partial}(A\otimes\mathbb{C}t^{-n}).$$ In particular, $\mathcal{L}_{(0)}=A\otimes \mathbb{C}t^{-1}+B/\partial A$. 

\vspace{0.2cm}

\noindent Let $\rho:L(A\oplus B)\rightarrow\mathcal{L}$ be a natural linear map defined by $$\rho( u\otimes t^n)=u\otimes t^n+\hat{\partial}L(A).$$ For $u\in A\oplus B$, $n\in\mathbb{Z}$, we set $u(n)=\rho(u\otimes t^n)$ and $u(z)=\sum_{n\in\mathbb{Z}}u(n)z^{-n-1}$. Let $W$ be a $\mathcal{L}$-module. We use $u_W(n)$ or sometimes just $u(n)$ for the corresponding operator on $W$ and we write $u_W(z)=\sum_{n\in\mathbb{Z}}u(n)z^{-n-1}\in (\End W)[[z,z^{-1}]]$. The commutator relations in terms of generating functions are the following:

\begin{eqnarray}
&&[a(z_1),a'(z_2)]=0,\label{vaaa'}\\
&&[a(z_1), b(z_2)]=z_2^{-1}\delta\left(\frac{z_1}{z_2}\right)(a_0b)(z_2),\label{vaab}\\
&&[b(z_1),b'(z_2)]=z_2^{-1}\delta\left(\frac{z_1}{z_2}\right)(b_0b')(z_2)+(b_1b')(z_2)\frac{\partial}{\partial z_2}z_2^{-1}\delta\left(\frac{z_1}{z_2}\right)\label{vabb'}
\end{eqnarray}
for $a,a'\in A$, $b,b'\in B$. 
\vspace{0.2cm}

\noindent Next, we define $\mathcal{L}^{\geq 0}=\rho((A\oplus B)\otimes \mathbb{C}[t])\subset \mathcal{L}$, and $\mathcal{L}^{<0}=\rho((A\oplus B)\otimes t^{-1}\mathbb{C}[t^{-1}])\subset \mathcal{L}$. As a vector space, $\mathcal{L}=\mathcal{L}^{\geq 0}\oplus \mathcal{L}^{<0}$. The subspaces $\mathcal{L}^{\geq 0}$ and $\mathcal{L}^{<0}$ are graded sub-algebras.

\vspace{0.2cm}

\noindent We now consider $\mathbb{C}$ as the trivial $\mathcal{L}^{\geq 0}$-module and form the following induced module $$V_{\mathcal{L}}=U(\mathcal{L})\otimes_{U(\mathcal{L}^{\geq 0})}\mathbb{C}.$$ 
In view of the Poincare-Birkhoff-Witt theorem, we have $V_{\mathcal{L}}=U(\mathcal{L}^{<0})$ as a vector space. We set ${\bf 1}=1\in V_{\mathcal{L}}$. We may consider $A\oplus B$ as a subspace: $$A\oplus B\rightarrow V_{\mathcal{L}}, \ \  a+b\mapsto a(-1){\bf 1}+b(-1){\bf 1}.$$ 
We assign $deg~\mathbb{C}=0$. Then $V_{\mathcal{L}}=\oplus_{n\in\mathbb{N}}(V_{\mathcal{L}})_{(n)}$ is a restricted $\mathbb{N}$-graded $\mathcal{L}$-module. 
\begin{prop} \cite{FKRW, MeP} 

\vspace{0.2cm}

\noindent There exists a unique vertex algebra structure on $V_{\mathcal{L}}$ with $Y(u,x)=u(x)$ for $u\in A\oplus B$. In fact, the vertex algebra $V_{\mathcal{L}}$ is a $\mathbb{N}$-graded vertex algebra and it is generated by $A\oplus B$. Furthermore, any restricted $\mathcal{L}$-module $W$ is naturally a $V_{\mathcal{L}}$-module with $Y_W(u,x)=u_W(x)$ for $u\in A\oplus B$. Conversely, any $V_{\mathcal{L}}$-module $W$ is naturally a restricted $\mathcal{L}$-module with $u_W(x)=Y_W(u,x)$ for $u\in A\oplus B$.
\end{prop}

\vspace{0.2cm} 

\noindent Now, we set 
\begin{eqnarray*}
&&E_0=Span\{\mathfrak{e}-{\bf 1},a(-1)a'-a*a'~|~a,a'\in A\}\subset (V_{\mathcal{L}})_{(0)},\\
&&E_1=Span\{a(-1)b-a\cdot b~|~a\in A,b\in\mathfrak{L}\}\subset (V_{\mathcal{L}})_{(1)},\\
&&E=E_0\oplus E_1.
\end{eqnarray*}
We define 
$$I_{B}=U(\mathcal{L})\mathbb{C}[D]E.$$ The vector space $I_{B}$ is an $\mathcal{L}$-submodule of $V_{\mathcal{L}}$. We set $$V_{B}=V_{\mathcal{L}}/I_{B}.$$

\begin{prop}\cite{GMS, LiY} \ \ 

\vspace{0.2cm}

\noindent (i) $V_{B}$ is a $\mathbb{N}$-graded vertex algebra such that $(V_{B})_{(0)}=A$ and $(V_{B})_{(1)}=B$ (under the linear map $v\mapsto v(-1){\bf 1}$) and $V_{B}$ as a vertex algebra is generated by $A\oplus B$. Furthermore, for any $n\geq 1$,
\begin{eqnarray*}
&&(V_{B})_{(n)}\\
&&=span\{b_1(-n_1).....b_k(-n_k){\bf 1}~|~b_i\in B,n_1\geq...\geq n_k\geq 1, n_1+...+n_k=n\}.
\end{eqnarray*}

\vspace{0.2cm}

\noindent (ii) A $V_{B}$-module $W$ is a restricted module for the Lie algebra $\mathcal{L}$ with $v(n)$ acting as $v_n$ for $v\in A\oplus B$, $n\in\mathbb{Z}$. Furthermore, the set of $V_{B}$-submodules is precisely the set of $\mathcal{L}$-submodules.
\end{prop}

\vspace{0.2cm}

\noindent Next, we recall definition of Lie algebroid and its module. Also, we will review construction of $V_B$-modules from modules of Lie $A$-algebroid $B/A\partial(A)$. 

\begin{dfn} Let $A$ be a commutative associative algebra. A {\em Lie $A$-algebroid} is a Lie algebra $\mathfrak{g}$ equipped with an $A$-module structure and a module action of $\mathfrak{g}$ on $A$ by derivation such that 
$$[u,av]=a[u,v]+(ua)v,\ \ a(ub)=(au)b$$
for all $u,v\in \mathfrak{g},a,b\in A$.

\vspace{0.2cm}

\noindent A {\em module for a Lie $A$-algebroid} $\mathfrak{g}$ is a vector space $W$ equipped with a $\mathfrak{g}$-module structure and an $A$-module structure such that 
$$u(aw)-a(uw)=(ua)w,~a(uw)=(au)w$$ for $a\in A$, $u\in\mathfrak{g}$, $w\in W$. 
\end{dfn}

\begin{prop}\label{prop30}\cite{LiY} 
Let $W=\oplus_{n\in\mathbb{N}}W_{(n)}$ be a $\mathbb{N}$-graded $V_{B}$-module with $W_{(0)}\neq\{ 0\}$. Then $W_{(0)}$ is an $A$-module with $a\cdot w=a_{-1}w$ for $a\in A$, $w\in W_{(0)}$ and $W_{(0)}$ is a module for the Lie algebra $B/A\partial(A)$ with $b\cdot w=b_0w$ for $b\in B$, $w\in W_{(0)}$. Furthermore, $W_{(0)}$ equipped with these module structures is a module for Lie $A$-algebroid $B/A\partial A$. If $W$ is graded simple, then $W_{(0)}$ is a simple module for Lie $A$-algebroid $B/A\partial A$.

\end{prop}

\vspace{0.2cm} 

\noindent Now, we set $\mathcal{L}_{\pm }=\oplus_{n\geq 1}\mathcal{L}_{(\pm n)}$ and $\mathcal{L}_{\leq 0}=\mathcal{L}_{-}\oplus \mathcal{L}_{(0)}$. Let $U$ be a module for the Lie algebra $\mathcal{L}_{(0)}(=A\oplus B/\partial (A))$. Then $U$ is an $\mathcal{L}_{(\leq 0)}$-module under the following actions:
$$a(n-1)\cdot u=\delta_{n,0}au,\ \ b(n)\cdot u=\delta_{n,0}u\text{ for }a\in A,b\in B, n\geq 0.$$ 
Next, we form the induced $\mathcal{L}$-module $M(U)=\Ind_{\mathcal{L}_{(\leq 0)}}^{\mathcal{L}}U$. Endow $U$ with degree 0, making $M(U)$ a $\mathbb{N}$-graded $\mathcal{L}$-module. In fact, $M(U)$ is a $V_{\mathcal{L}}$-module. 
We set $$W(U)=span\{v_nu~|~v\in E, ~n\in\mathbb{Z}, ~u\in U\}\subset M(U),$$ and
$$M_{B}(U)=M(U)/U(\mathcal{L})W(U).$$
\begin{prop}\label{simplemodulerelations}\cite{LiY} \ \ 

\vspace{0.2cm} 

\noindent (i) Let $U$ be a module for the Lie algebra $\mathcal{L}_{(0)}$. Then $M_{B}(U)$ is a $V_{B}$-module. If $U$ is a module for the Lie $A$-algebroid $B/A\partial A$ then $(M_{B}(U))_{(0)}=U$.

\vspace{0.2cm} 

\noindent (ii) Let $U$ be a module for the Lie $A$-algebroid $B/A\partial A$. Then there exists a unique maximal graded $U(\mathcal{L})$-submodule $J(U)$  of $M(U)$ with the property that $J(U)\cap U=0$. Moreover, $L(U)=M(U)/J(U)$ is a $\mathbb{N}$-graded $V_{B}$-module such that $L(U)_{(0)}=U$ as a module for the Lie $A$-algebroid $B/A\partial A$. If $U$ is a simple $B/A\partial A$, $L(U)$ is a graded simple $V_{B}$-module.

\vspace{0.2cm}

\noindent (iii) Let $W=\coprod_{n\in\mathbb{N}}W_{(n)}$ be an $\mathbb{N}$-graded simple $V_B$-module with $W_{(0)}\neq 0$. Then $W\cong L(W_{(0)})$. 

\vspace{0.2cm}

\noindent (iv) For any complete set $H$ of representatives of equivalence classes of simple modules for the Lie $A$-algebroid $B/A\partial A$, $\{L(U)~|~U\in H\}$ is a complete set of representatives of equivalence classes of simple $\mathbb{N}$-graded simple $V_B$-modules.
\end{prop}

\subsection{Indecomposable $C_2$-cofinite $\mathbb{N}$-graded vertex algebras associated with simple Leibniz algebras}

\ \ \  

\vspace{0.2cm}

\noindent In this subsection, we construct an indecomposable $C_2$-cofinite $\mathbb{N}$-graded vertex algebras from simple Leibniz algebras. Also, we study irreducible modules of these vertex algebras and examine a connection between these irreducible modules and irreducible modulese of rational affine vertex operator algebras. 

\vspace{0.1cm}

\noindent Let $(A,*)$ be a finite dimensional commutative associative algebra with the identity $\hat{1}$, and let $B$ be a finite dimensional vertex $A$-algebroid such that 

\vspace{0.1cm}

\noindent (i) $B$ is a simple Leibniz algebra that is not a Lie algebra. Consequently, there exists a simple Lie algebra $\mathfrak{g}$ such that $B=\mathfrak{g}+\partial(A)$ and $\partial(A)$ is a simple $\mathfrak{g}$-module;

\vspace{0.1cm}

\noindent (ii) there exist $u,v\in\mathfrak{g}$ such that $u_1v\neq 0$, $Ker(\partial)=\mathbb{C}\hat{1}$ and $N\cdot B\subset \partial(A)$. 

\vspace{0.1cm}

\noindent By Theorem \ref{vertexalgebroidsimpleLeibnizalgebra}, there exists a dominant integral $\lambda$ such that $B\cong B_{\lambda}$ as vertex algbroids. For simplicity, from now on we will assume that $$A:=A_{\lambda}=\mathbb{C}\hat{1}+N_{\lambda},~B:=B_{\lambda}=\mathfrak{g}+\partial(N_{\lambda}).$$ Here $\mathfrak{g}$ is a Lie algebra of $ADE$-type that is not $E_8$. The vertex $A_{\lambda}$-algebroid $B_{\lambda}$ satisfies the following properties: for $g,g'\in\mathfrak{g}$, $a,a'\in N_{\lambda}$, we have 
\begin{eqnarray*}
&&g_1g'=\langle g,g'\rangle\hat{1},~a*a'=0,~a\cdot (g+\partial(a'))=\partial(g_0a),\\
&&g_1\partial(a)=g\cdot a=\partial(a)_1g,~\partial(a)_1\partial(a')=0.
\end{eqnarray*} Here, $\langle ~,~\rangle$ is the Killing form of $\mathfrak{g}$. Following subsection 4.1, we let $\mathfrak{h}$ be a Cartan subalgebra of $\mathfrak{g}$ and identify $\mathfrak{h}$ with $\mathfrak{h}^*$ by means of the Killing form $\langle~,~\rangle$. Let $\Delta$ be the root system of $\mathfrak{g}$ with respect to $\mathfrak{h}$. We have the root space decomposition $\mathfrak{g}=\mathfrak{h}\oplus\sum_{\alpha\in\Delta}\mathfrak{g}_{\alpha}$. In addition, for each $\alpha\in \Delta$, we set $\mathfrak{g}^{\alpha}=\mathfrak{g}_{\alpha}\oplus\mathbb{C}h_{\alpha}\oplus\mathfrak{g}_{-\alpha}$ and denote the set of positive roots by $\Delta_+$. Also, we let $\theta$ be the highest root. Hence $\langle\theta,\theta\rangle=2$. We will fix vectors $e_{\theta}\in\mathfrak{g}_{\theta}$, $f_{\theta}\in\mathfrak{g}_{-\theta}$ such that $(e_{\theta})_1f_{\theta}=\langle e_{\theta},f_{\theta}\rangle\hat{1}=\hat{1}$. So, we have 
$$(h_{\theta})_0e_{\theta}=2e_{\theta},~(h_{\theta})_0f_{\theta}=-2f_{\theta},~(e_{\theta})_0f_{\theta}=h_{\theta}.$$ By Lemma \ref{dominantintegral}, the following statements hold:
 
\vspace{0.1cm}

\noindent (i) For $\alpha\in \Delta_+$, $\lambda(h_{\alpha})=1$. In particular $\lambda(h_{\theta})=1$.

\vspace{0.1cm}

\noindent (ii) We let $v_{\lambda}\in N_{\lambda}$ such that $h\cdot v_{\lambda}=\lambda(h)v_{\lambda}$, $g\cdot v_{\lambda}=0$ for all $h\in \mathfrak{h},g\in\mathfrak{g}_+$. For $\alpha\in \Delta_+$, $a\in N_{\lambda}$, we have $$(f_{\alpha})_0(f_{\alpha})_0a=0\text{ and }(e_{\alpha})_0(e_{\alpha})_0a=0.$$

\vspace{0.2cm}

\noindent Now, we let $V_{B_{\lambda}}$ be the $\mathbb{N}$-graded vertex algebra associated with the vertex $A_{\lambda}$-algebroid $B_{\lambda}$. We will use the following Proposition to help us show that $V_{B_{\lambda}}$ is an indecomposable non-simple vertex algebra. 

\begin{prop}\label{indecom}\cite{JY} Let $V=\oplus_{n=0}^{\infty}V_{(n)}$ be a $\mathbb{N}$-graded vertex algebra that satisfies the following properties:

\noindent (a) $2\leq \dim V_{(0)}<\infty$, $1\leq \dim V_{(1)}<\infty$, $V$ is generated by $V_{(0)}$ and $V_{(1)}$;
\noindent (b) $V_{(0)}$ is a local algebra. 

\noindent Assume that one of the following statements hold.

\noindent (i) An ideal generated by $rad\langle~,~\rangle$ is not zero;

\noindent (ii) $V_{(0)}$ is not a simple module for a Lie $V_{(0)}$-algebroid $V_{(1)}/Ann_{V_{(1)}}(V_{(0)})$;

\noindent (iii) $V_{(0)}$ is not a simple module for a Lie $V_{(0)}$-algebroid $V_{(1)}/V_{(0)}D(V_{(0)})$.

\noindent Then $V$ is an indecomposable non-simple vertex algebra.
\end{prop}

\noindent Now, we will investigate the $\mathbb{N}$-graded vertex algebra $V_{B_{\lambda}}$ by using Proposition \ref{indecom}. Notice that when $V=V_{B_{\lambda}}$, we have 

\noindent 1) $V_{(1)}=\mathfrak{g}+\partial(N_{\lambda})$, $V_{(0)}D(V_{(0)})=\partial(N_{\lambda})$, and $V_{(0)}=A_{\lambda}(=\mathbb{C}\hat{1}+N_{\lambda})$. 

\noindent 2) Moreover, $V_{(0)}=A_{\lambda}$ is a local algebra, $V_{(1)}/V_{(0)}D(V_{(0)})=(\mathfrak{g}+\partial(N_{\lambda})/\partial(N_{\lambda})\cong\mathfrak{g}$ as Lie algebras, $N_{\lambda}$ is a module for the Lie algebra $(\mathfrak{g}+\partial(N_{\lambda})/\partial(N_{\lambda})$ and $N_{\lambda}$ is an ideal of the commutative associative algebra $A_{\lambda}$.

\vspace{0.1cm} 

\noindent Now, let $u+\partial(A_{\lambda})\in(\mathfrak{g}+\partial(N_{\lambda})/\partial(N_{\lambda})$, and let $\beta^1 \hat{1}+a'\in A_{\lambda}$. Here $a'\in N_{\lambda}$, $\beta\in\mathbb{C}$ . For $a''\in N_{\lambda}$, we have
\begin{eqnarray*}
&&(u+\partial(A_{\lambda}))_0((\beta\hat{1}+a')\cdot a'')-(\beta\hat{1}+a')\cdot ((u+\partial(A_{\lambda}))_0a'')\\
&&=(u+\partial(A_{\lambda}))_0\beta a''-(\beta\hat{1}+a')\cdot(u_0a'')\\
&&=u_0(\beta a'')-\beta u_0a''\\
&&=0\\
&&=(u_0a')*a''\\
&&=((u+\partial(A_{\lambda}))_0(\beta\hat{1}+a'))\cdot a''
\end{eqnarray*}
and 
\begin{eqnarray*}
&&(\beta\hat{1}+a'')\cdot ((u+\partial(A_{\lambda}))_0a'')\\
&&=(\beta\hat{1}+a'')\cdot  u_0a''\\
&&=\beta u_0a''\\
&&=(\beta u+\partial(A_{\lambda}))_0a''\\
&&=(( \beta\hat{1}+a'')\cdot u+\partial(A_{\lambda}))_0 a''.
\end{eqnarray*}
Hence, we can conclude that $N_{\lambda}$ is a module for the Lie $A_{\lambda}$-algebroid $(\mathfrak{g}+\partial(N_{\lambda})/\partial(N_{\lambda})$. Therefore, $A_{\lambda}$ is not a simple module for the Lie $A_{\lambda}$-algebroid $(\mathfrak{g}+\partial(N_{\lambda})/\partial(N_{\lambda})$. By Proposition \ref{indecom}, we can conclude that $V_{B_{\lambda}}$ is an indecomposable non-simple $\mathbb{N}$-graded vertex algebra. We summarize this discussion in the following theorem.

\vspace{0.2cm} 

\begin{thm} Let $\mathfrak{g}$ be a simple Lie algebra of $ADE$-type that is not $E_8$. Let $\langle~,~\rangle$ be the killing form of $\mathfrak{g}$ and let $\mathfrak{h}$ be a Cartan subalgebra of $\mathfrak{g}$. For a dominant integral $\lambda\in \mathfrak{h}^*$, we let $A_{\lambda}=\mathbb{C}\hat{1}\oplus N_{\lambda}$ be a vector space such that $N_{\lambda}$ is an irreducible $\mathfrak{g}$-module. We set $B_{\lambda}=\mathfrak{g}+\partial(A_{\lambda}).$ Assume that $A_{\lambda}$ is a unital commutative associative algebra and $B_{\lambda}$ is a vertex $A_{\lambda}$-algebroid that satisfy Proposition \ref{simpleLiealgebraVAB}. Then the $\mathbb{N}$-graded vertex algebra $V_{B_{\lambda}}$ is an indecomposable nonsimple vertex algebra. 
\end{thm}

\begin{lem}\label{simplemoduleLiealgebroid} Let $U$ be an irreducible module for $B_{\lambda}/A_{\lambda}\partial(A_{\lambda})$ viewed as a Lie algebra. Naturally view $U$ a module for $B_{\lambda}$ as a simple Leibniz algebra. Define an action of $A_{\lambda}$ on $U$ by letting  $\hat{1}$ acts as a scalar 1 and elements in $N_{\lambda}$ act trivially on $U$. Then $U$ is an irreducible module for the Lie $A_{\lambda}$-algbroid $B_{\lambda}/A_{\lambda}\partial(A_{\lambda})$.
\end{lem}
\begin{proof} Let $U$ be an irreducible $B_{\lambda}/A_{\lambda}\partial(A_{\lambda})$-module. To show that $U$ is an irreducible module for a Lie $A_{\lambda}$-algbroid $B_{\lambda}/A_{\lambda}\partial(A_{\lambda})$, we will follow the proof of Lemma 33 in \cite{JY2}. Recall that there exists a unique algebra homomorphism $\rho:A_{\lambda}\rightarrow\mathbb{C}$ with $Ker(\rho)$ being the unique maximal ideal $N_{\lambda}$ because $A_{\lambda}$ is a local algebra. Hence, any vector space can be $A_{\lambda}$-module. In particular, $U$ is a module of $A_{\lambda}$. Now, to show that $U$ is a module for the Lie $A_{\lambda}$-algebroid $B_{\lambda}/A_{\lambda}\partial(A_{\lambda})$, we need to show that for $a\in A_{\lambda}, b\in B_{\lambda}/A_{\lambda}\partial(A_{\lambda}), w\in U$, $$b_0(a\cdot w)-a\cdot (b_0w)=(b_0a)\cdot w\text{ and }a\cdot (b_0w)=(a\cdot b)_0w.$$ Let $\beta\hat{1}+a'\in A_{\lambda}$, $u+\partial(A_{\lambda})\in B_{\lambda}/A_{\lambda}\partial(A_{\lambda}), w\in U$. Here, $\beta\in\mathbb{C}$, $a'\in N_{\lambda}$. Since
\begin{eqnarray*}
&&(u+\partial(A_{\lambda}))_0((\beta\hat{1}+a')\cdot w)-(\beta\hat{1}+a')\cdot ((u+\partial(A_{\lambda}))_0w)\\
&&=(u+\partial(A_{\lambda}))_0(\beta w)-(\beta\hat{1}+a')\cdot (u_0w)\\
&&=\beta u_0w-\beta u_0w \\
&&=0\\
&&=(u_0a')\cdot w\\
&&=((u+\partial(A_{\lambda}))_0(\beta\hat{1}+a'))\cdot w
\end{eqnarray*}
and
\begin{eqnarray*}
&&(\beta\hat{1}+a')\cdot ( (u+\partial(A_{\lambda}))_0w)\\
&&=(\beta\hat{1}+a')\cdot u_0w\\
&&=\beta u_0w\\
&&=(\beta u+\partial(A_{\lambda}))_0w\\
&&=((\beta\hat{1}+a')\cdot (u+\partial(A_{\lambda})))_0w,
\end{eqnarray*}
we can conclude that $U$ is a module for the Lie $A_{\lambda}$-algebroid $B_{\lambda}/A_{\lambda}\partial(A_{\lambda})$. It is clear that $U$ is an irreducible module for the Lie $A_{\lambda}$-algebroid $B_{\lambda}/A_{\lambda}\partial(A_{\lambda})$.
\end{proof}

\begin{lem}\label{1-1Liever} There is a one-to-one correspondence between the set of representatives of equivalence classes of simple $\mathfrak{g}$-modules and the set of representatives of equivalence classes of $\mathbb{N}$-graded simple $V_{B_{\lambda}}$-modules.
\end{lem}
\begin{proof} By following the proof of Lemma 34 in \cite{JY2}, we can prove that if $W$ is an irreducible module for the Lie $A_{\lambda}$-algebroid $B_{\lambda}/A_{\lambda}\partial(A_{\lambda})$, then $N_{\lambda} W=Span\{a\cdot w~|~a\in N_{\lambda}, w\in W\}=\{0\}$ and $W$ is a irreducible $\mathfrak{g}$-module. By using Lemma \ref{simplemoduleLiealgebroid} and following the proof of Lemma 35 in \cite{JY2}, we can show that there is a one-to-one correspondence between the set of representatives of equivalence classes of simple $\mathfrak{g}$-modules and the set of representatives of equivalence classes of $\mathbb{N}$-graded simple $V_{B_{\lambda}}$-modules.
\end{proof}
\noindent Now, we let $(e_{\theta}(-1)e_{\theta})$ be an ideal of $V_{B_{\lambda}}$ that is generated by $e_{\theta}(-1)e_{\theta}$.

\begin{lem}\label{idealetheta} $(e_{\theta}(-1)e_{\theta})\cap A_{\lambda}=\{0\}$ and $(e_{\theta}(-1)e_{\theta})\cap B_{\lambda}=\{0\}$.
\end{lem}
\begin{proof} The proof of this lemma is very similar to the proof of Lemma 36 in \cite{JY2}. In fact, by substituting $e$ in Lemma 36 in \cite{JY2} by $e_{\theta}$, and using Theorem \ref{sl2insideB} and the fact that $e_{\theta},f_{\theta},h_{\theta}\in B_{\lambda}$, one can show that $a(n)e_{\theta}(-1)e_{\theta}=0$ for all $a\in A_{\lambda}$. 

\vspace{0.2cm}

\noindent Now, we will show that for $b \in B_{\lambda}$, $b(n)e_{\theta}(-1)e_{\theta}=0$ for all $n\geq 1$. Note that by Lemma \ref{dominantintegral} we have $(e_{\theta})_0(e_{\theta})_0a'=0$ for all $a'\in A_{\lambda}$. Let $b\in B_{\lambda}$. Then there exist $\beta_f,\beta_h\in\mathbb{C}, g\in \mathfrak{g}~\backslash~ (\mathbb{C}f_{\theta}\oplus \mathbb{C}h_{\theta})$ and $a\in A_{\lambda}$ such that $b=\beta_ff_{\theta}+\beta_h h_{\theta}+g+\partial(a)$. Since 
\begin{eqnarray*}
b_0e_{\theta}&&=\beta_f(-h_{\theta})+2\beta_h e_{\theta}+g_0e_{\theta}\\
b_1e_{\theta}&&=\beta_f{\bf 1} +(e_{\theta})_0a\\
\end{eqnarray*}
we have 
\begin{eqnarray*}
b(1)e_{\theta}(-1)e_{\theta}&&=e_{\theta}(-1)b(1)e_{\theta}+(b_0e_{\theta})(0)e_{\theta}+(b_1e_{\theta})(-1)e_{\theta}\\
&&=e_{\theta}(-1)(\beta_f {\bf 1}+(e_{\theta})_0a)+(\beta_f(-h_{\theta})+2\beta_h e_{\theta}+g_0e_{\theta} )(0)e_{\theta}\\
&&\ \ \  +(\beta_f{\bf 1} +(e_{\theta})_0a  )(-1)e_{\theta}\\
&&=2\beta_f e_{\theta}(-1){\bf 1}+((e_{\theta})_0a)(-1)e_{\theta}-\partial(((e_{\theta})_0a)_0e_{\theta})\\
&&\ \ \ \ -2\beta_fe_{\theta}+(g_0e_{\theta})_0e_{\theta}+\partial((e_{\theta})_0(e_{\theta})_0a)\\
&&=-(e_{\theta})_0g_0(e_{\theta})\\
&&=0
\end{eqnarray*}
and
\begin{eqnarray*}
b(2)(e_{\theta})(-1)e_{\theta}&&=e_{\theta}(-1)b(2)e+(b_0e_{\theta})(1)e_{\theta}+2(b_1e_{\theta})(0)e_{\theta}\\
&&=(\beta_f(-h_{\theta})+2\beta_h e_{\theta}+g_0e_{\theta})(1)e_{\theta}\\
&&=0.
\end{eqnarray*}
Consequently, we have $b(n)(e_{\theta})(-1)e_{\theta}=0$ for all $n\geq 1$.  Now, by substituting $e$ in Lemma 36 in \cite{JY2} by $e_{\theta}$ and following the proof of Lemma 36 in \cite{JY2} after equation (27) step by step, one can show that $(e_{\theta}(-1)e_{\theta})\cap A_{\lambda}=\{0\}$ and $(e_{\theta}(-1)e_{\theta})\cap B_{\lambda}=\{0\}$ as desired.
\end{proof}

\noindent Now, we set $$\overline{V_{B_{\lambda}}}=V_{B_{\lambda}}/(e_{\theta}(-1)e_{\theta}).$$ Clearly, by Proposition \ref{indecom} and Lemma \ref{idealetheta}, we can conclude immediately that 
\begin{prop}\label{indecompVB} $\overline{V_{B_{\lambda}}}=\oplus_{n=0}^{\infty}\left(\overline{V_{B_{\lambda}}}\right)_{(n)}$ is an indecomposable non-simple $\mathbb{N}$-graded vertex algebra such that $\left(\overline{V_{B_{\lambda}}}\right)_{(0)}=A_{\lambda}$ and $\left(\overline{V_{B_{\lambda}}}\right)_{(1)}=B_{\lambda}$. 
\end{prop}

\vspace{0.2cm}

\noindent Following \cite{JY2} page 809, we let $\hat{\mathfrak{g}}=\mathfrak{g}\otimes \mathbb{C}[t,t^{-1}]\oplus\mathbb{C}c$ be the affine Lie algebra where $c$ is central and 
$$[u\otimes t^m,v\otimes t^n]=[u,v]\otimes t^{m+n}+m(( u,v))\delta_{m+m,0}c.$$ 
Here $(( ~,~))$ is a symmetric invariant bilinear form of $\mathfrak{g}$ normalize so that the square length of long root is 2. The generalized Verma $\hat{\mathfrak{g}}$-module $M_{\mathfrak{g}}(k,0)$ is a vertex operator algebra. For $u\in \overline{V_{B_{\lambda}}}$, we set $Y_{\overline{V_{B_{\lambda}}}}(u,z)=\sum_{n\in\mathbb{Z}}u[n]z^{-n-1}$. Since $\mathfrak{g}$ is a subset of $\left(\overline{V_{B_{\lambda}}}\right)_{(1)}$ and $\mathfrak{g}$ is a Lie algebra with the killing form $\langle ~,~\rangle$ such that $\langle s,s'\rangle\hat{1}=s[1]s'$ for all $s,s'\in\mathfrak{g}$, the map $\hat{\mathfrak{g}}\rightarrow \End (\overline{V_B}): s\otimes t^m\mapsto s[m]$ is a representation of the affine Kac-Moody algebra $\hat{\mathfrak{g}}$ of level $k$ where $\langle s,s'\rangle=k((s,s'))$ for $s,s'\in\mathfrak{g}$. Since $\langle h_{\theta},h_{\theta}\rangle\hat{1}=h_{\theta}[1]h_{\theta}=2\hat{1}$ and $((h_{\theta},h_{\theta}))=2$, we then have that $k=1$. 

\vspace{0.2cm}

\noindent Let $U$ be the vertex sub-algebra of $\overline{V_{B_{\lambda}}}$ that is generated by $\mathfrak{g}$. By replacing $e$ on page 809 in \cite{JY2} by $e_{\theta}$ and following the proof on that page step by step, one can show that $U$ is isomorphic to the rational vertex operator algebra $L_{\hat{\mathfrak{g}}}(1,0)$ associated to the Lie algebra $\mathfrak{g}$. Moreover, $\overline{V_{B_{\lambda}}}$ is an integrable $\hat{\mathfrak{g}}$-module. Note that $L_{\hat{\mathfrak{g}}}(1,0)=M_{\mathfrak{g}}(1,0)/J_{1,0}$ where$J_{1,0}$ is the maximal proper submodule of $M_{\mathfrak{g}}(1,0)$ that is generated by the vector $e_{\theta}(-1)e_{\theta}(-1){\bf 1}$.

\vspace{0.2cm}

\noindent Recall that for a long root $\alpha\in \Delta$ and let $e\in \mathfrak{g}_{\alpha}$, $Y(e,z)^2=0$ on $L_{\hat{\mathfrak{g}}}(1,0)$. Moreover, for any module $W$ for $L_{\hat{\mathfrak{g}}}(1,0)$ viewed as a vertex algebra, $Y_W(e,z)^2=0$. In particular $Y_W(e_{\theta},z)^2=Y_W(f_{\theta},z)^2=0$ (cf. Proposition 6.6.5 of \cite{LLi}). Consequently, the following proposition holds.

\begin{prop} Let $\beta\in \Delta$ be a long root and let $e_{\beta}\in \mathfrak{g}_{\beta}$. Then $Y(e_{\beta},z)^2=0$ on $\overline{V_{B_{\lambda}}}$. Furthermore, for any module $W$ of $\overline{V_{B_{\lambda}}}$, $Y_W(e_{\beta},z)^2=0$. In particular, $Y_W(e_{\theta},z)^2=Y_W(f_{\theta},z)^2=0$.
\end{prop}

\begin{lem} $\overline{V_{B_{\lambda}}}$ satisfies the $C_2$-condition.
\end{lem}
\begin{proof} Observe that 
\begin{eqnarray*}
&&\overline{V_{B_{\lambda}}}/C_2(\overline{V_{B_{\lambda}}}  )\\
&&=Span\{a+C_2(\overline{V_{B_{\lambda}}}) , b+C_2(\overline{V_{B_{\lambda}}}) , b^1[-1]...b^k[-1]{\bf 1}+C_2(\overline{V_{B_{\lambda}}})~|~a\in A_{\lambda}, \\
&&\ \ \ \ \ \ \ \ \  \ \ \ \ \ \  b,b^i\in \mathfrak{g}, k\geq 2\}.
\end{eqnarray*} 
Next, by following the proof of Proposition 12.6 in \cite{DLM}, we can show that $\dim \overline{V_{B_{\lambda}}}/C_2(\overline{V_{B_{\lambda}}}  )<\infty$.
\end{proof}

\vspace{0.2cm}

\noindent Next, we will study $\mathbb{N}$-graded $\overline{V_{B_{\lambda}}}$-modules. Observe that $A_{\lambda}\oplus B_{\lambda}$ generates $\overline{V_{B_{\lambda}}}$ as a vertex algebra. Consequently, if $W$ is a $\overline{V_{B_{\lambda}}}$-module, then $W$ is a restricted $\mathcal{L}$-module with $u(n)$ acting as $u_n$ for $u\in A_{\lambda}\oplus B_{\lambda}$, $n\in\mathbb{Z}$. Moreover, the set of $\overline{V_{B_{\lambda}}}$-submodules is the set of $\mathcal{L}$-submodules. 

\vspace{0.2cm}

\noindent Now, we will recall the following results in \cite{LiY}. 
\begin{prop}\cite{LiY} Let $(W,Y_W)$ be a $V_{\mathcal{L}}$-module. Assume that for any $a,a'\in A, b\in B$,
\begin{eqnarray*}
&&Y_W(\mathfrak{e},z)u=u,\\
&&Y_W(a(-1)a',z)u=Y_W(a*a',z)u,\\ 
&&Y_W(a(-1)b,z)u=Y_W(a\cdot b,z)u,
\end{eqnarray*}
for all $u\in U$ where $U$ is a generating subspace of $W$ as a $V_{\mathcal{L}}$-module, then $W$ is naturally a ${V_{B_{\lambda}}}$-module. 
\end{prop}

\begin{prop}\label{annihilateVW}\cite{LiY} Let $V$ be a vertex algebra and let $I$ be a (two-sided) ideal generated by a subset $S$. Let $(W,Y_W)$ be a $V$-module  and let $U$ be a generating subspace of $W$ as a $V$-module such that $Y_W(v,x)u=0$ for $v\in S$, $u\in U$. Then $Y_W(v,x)=0$ for $v\in I$.
\end{prop}
\begin{lem} Let $(Q,Y_Q)$ be a $V_{B_{\lambda}}$-module such that $Y(e_{\theta}(-1)e_{\theta},z)u=0$ for all $u\in F$ where $F$ is a generating subspace of $Q$ as a $V_{B_{\lambda}}$-module. Then $Q$ is a $\overline{V_{B_{\lambda}}}$-module. 
\end{lem}
\begin{proof} By using Proposition \ref{annihilateVW}, one can obtain the above statement very easily. 
\end{proof}

\begin{lem} Let $W=\oplus_{n=0}^{\infty}W_{(n)}$ be a $\mathbb{N}$-graded $\overline{V_{B_{\lambda}}}$-module with $W_{(0)}\neq \{0\}$. Then 

\vspace{0.1cm}

\noindent (i) $W_{(0)}$ is an $A_{\lambda}$-module with $a\cdot w=a_{-1}w$ for $a\in A_{\lambda}$, $w\in W_{(0)}$ and $W_{(0)}$ is a module for the Lie algebra $B_{\lambda}/A_{\lambda}\partial(A_{\lambda})(\cong \mathfrak{g})$ with $b\cdot w=b_0w$ for $b\in B_{\lambda}$, $w\in W_{(0)}$. Furthermore, $W_{(0)}$ equipped with these module structure is a module for the Lie $A_{\lambda}$-algebroid $B_{\lambda}/A_{\lambda}\partial(A_{\lambda})$. 

\noindent (ii) Let $\beta\in\Delta$ be a long root and $e\in \mathfrak{g}_{\beta}$. Then 
$$e_0(e_0w)=0,~e_{-1}(e_{-1}w)=0.$$ In particular, $(e_{\theta})_0(e_{\theta})_0w=0$, $(e_{\theta})_{-1}((e_{\theta})_{-1}w)=0$, $(f_{\theta})_0((f_{\theta})_0w)=0$, and $(f_{\theta})_{-1}((f_{\theta})_{-1}w)=0$.

\noindent (iii) If $W$ is simple then $W_{(0)}$ is a finite dimensional irreducible module for the Lie $A_{\lambda}$-algebroid $B_{\lambda}/A_{\lambda}\partial(A_{\lambda})$. For simplicity, we assume that $W_{(0)}$ equals $L(\gamma)$ where $L(\gamma)$ is a finite dimensional irreducible module for the Lie algebra $B_{\lambda}/A_{\lambda}\partial(A_{\lambda})$. Then $\gamma$ is dominant integral and $\gamma(h_{\theta})\leq 1$.  
\end{lem}
\begin{proof} Following the proof of statements (i), (ii) in Lemma 42 of \cite{JY2} step by step, one can easily obtain the statments (i), (ii) of this Lemma.  Next, following the proof of Proposition 4.8 of \cite{LiY}, one can show that if $W$ is simple, then $W_{(0)}$ is an irreducible module for the algebra $B_{\lambda}/A_{\lambda}\partial(A_{\lambda})$. 

\vspace{0.1cm}

\noindent Now, we recall Lemma 6.6.8 of \cite{LLi}: let $U$ be a $\mathfrak{g}$-module on which $\mathfrak{h}$ acts semisimply. Suppose that $\mathfrak{g}_{\alpha}$,and $\mathfrak{g}_{-\alpha}$ act nilpotently on $U$ for some $\alpha\in\Delta$. Then for each $\beta\in \Delta$, $\mathfrak{g}_{\beta}$ acts nilpotently on $U$ and $U$ is a direct sum of finite-dimensional irreducible $\mathfrak{g}$-modules. 

\vspace{0.1cm}

\noindent By the above Lemma, we can conclude that $W(0)$ is finite dimensional. Let $w_0$ be a highest weight vector of $W(0)$ of uniquely determined weight $\gamma$. The submodule that is generated by $w_0$ must be all of $W(0)$. We set $\mathfrak{g}^{\theta}=\mathfrak{g}_{\theta}\oplus\mathbb{C}h_{\theta}\oplus \mathfrak{g}_{-\theta}$. $U(\mathfrak{g}^{\theta})w_0$ is an irreducible module for the three-dimensional simple Lie algebra $\mathfrak{g}^{\theta}$ with highest weight $\lambda(h_{\theta})$, so that $\dim U(\mathfrak{g}^{\theta})w_0=\gamma(h_{\theta})+1$. Since $(e_{\theta})_0(e_{\theta})_0 W_{(0)}=0$, we must have that $\gamma(h_{\theta})+1\leq 2$. Hence, $\gamma(h_{\theta})\leq 1$. 

\end{proof}
\begin{lem}\label{classifyirreduciblemodules}  \ \ 

\vspace{0.1cm}

\noindent (i) There exists a vertex algebra homomorphism $\tilde{\theta}:\overline{V_{B_{\lambda}}}\rightarrow  L_{\hat{\mathfrak{g}}}(1,0)$ such that $\tilde{\theta}(a)=\theta({a}){\bf 1}$ and $\tilde{\theta}(b)=b+A_{\lambda}\partial(A_{\lambda})$ for all $a\in A_{\lambda}$, $b\in B_{\lambda}$. 

\noindent (ii) If $W$ is an irreducible $L_{\hat{\mathfrak{g}}}(1,0)$-module then $W$ is an irreducible $\overline{V_{B_{\lambda}}}$-module.

\vspace{0.1cm}

\noindent (iii) If $U$ is an irreducible $\overline{V_{B_{\lambda}}}$-module then  $W$ is an irreducible $L_{\hat{\mathfrak{g}}}(1,0)$-module.
\end{lem}
\begin{proof} By following the proof of Lemma 44 in \cite{JY2}, we will obtain statements (i), (ii). By using Lemma \ref{1-1Liever} and statement (ii),we then have statement (iii). 
\end{proof}

\vspace{0.2cm}

\noindent We summarize results from Proposition \ref{indecompVB} to Lemma \ref{classifyirreduciblemodules} in the following theorem.
\begin{thm} Let $(A,*)$ be a finite dimensional commutative associative algebra with the identity $\hat{1}$, and let $B$ be a finite dimensional vertex $A$-algebroid such that 

\vspace{0.1cm}

\noindent (i) $B$ is a simple Leibniz algebra that is not a Lie algebra. Consequently, there exists a simple Lie algebra $\mathfrak{g}$ such that $B=\mathfrak{g}+\partial(A)$ and $\partial(A)$ is a simple $\mathfrak{g}$-module;

\vspace{0.1cm}

\noindent (ii) there exist $u,v\in\mathfrak{g}$ such that $u_1v\neq 0$, $Ker(\partial)=\mathbb{C}\hat{1}$ and $N\cdot B\subset \partial(A)$. 

\vspace{0.1cm}

\noindent Hence, there exists a dominant integral $\lambda$ such that $B\cong B_{\lambda}$ as vertex algbroid. For simplicity, from now on we will assume that $$A:=A_{\lambda}=\mathbb{C}\hat{1}+N_{\lambda},~B:=B_{\lambda}=\mathfrak{g}+\partial(N_{\lambda}).$$ Here $\mathfrak{g}$ is a Lie algebra of $ADE$-type that is not $E_8$.  

\vspace{0.2cm}

\noindent We will denote the Killing form of $\mathfrak{g}$ by $\langle ~,~\rangle$. We let $\mathfrak{h}$ be a Cartan subalgebra of $\mathfrak{g}$ and identify $\mathfrak{h}$ with $\mathfrak{h}^*$ by means of the Killing form $\langle~,~\rangle$. Let $\Delta$ be the root system of $\mathfrak{g}$ with respect to $\mathfrak{h}$. We have the root space decomposition $\mathfrak{g}=\mathfrak{h}\oplus\sum_{\alpha\in\Delta}\mathfrak{g}_{\alpha}$. We let $\theta$ be the highest root. Hence $\langle\theta,\theta\rangle=2$. We will fix vectors $e_{\theta}\in\mathfrak{g}_{\theta}$, $f_{\theta}\in\mathfrak{g}_{-\theta}$ such that $(e_{\theta})_1f_{\theta}=\langle e_{\theta},f_{\theta}\rangle\hat{1}=\hat{1}$. 

\vspace{0.2cm}

\noindent We let $V_{B_{\lambda}}$ be the $\mathbb{N}$-graded vertex algebra associated with the vertex $A_{\lambda}$-algebroid $B_{\lambda}$. We let $(e_{\theta}(-1)e_{\theta})$ be the ideal of $V_{B_{\lambda}}$ that is generated by $e_{\theta}(-1)e_{\theta}$. Also, we set $\overline{V_{B_{\lambda}}}=V_{B_{\lambda}}/(e_{\theta}(-1)e_{\theta}) $. The following statements hold:

\vspace{0.2cm}

\noindent i) $\overline{V_{B_{\lambda}}}=\oplus_{n=0}^{\infty}\left(\overline{V_{B_{\lambda}}}\right)_{(n)}$ is a $C_2$-cofinite indecomposable non-simple $\mathbb{N}$-graded vertex algebra such that $\left(\overline{V_{B_{\lambda}}}\right)_{(0)}=A_{\lambda}$ and $\left(\overline{V_{B_{\lambda}}}\right)_{(1)}=B_{\lambda}$. 

\noindent ii) Let $L_{\hat{\mathfrak{g}}}(1,0)$ be the rational affine vertex operator algebra associated with the simple Lie algebra $\mathfrak{g}$. Then a vector space $W$ is an irreducible $\mathbb{N}$-graded $\overline{V_{B_{\lambda}}}$-module if and only if $W$ is an irreducible $L_{\hat{\mathfrak{g}}}(1,0)$-module. 
 
\end{thm}

\subsection{Conformal vectors of $\overline{V_{B_{\lambda}}}$}

\ \ 

\vspace{0.2cm}

\noindent In this subsection, we use Virasoro elements of $L_{\hat{\mathfrak{g}}}(1,0)$ to study conformal vectors of $\overline{V_{B_{\lambda}}}$.

\begin{dfn}\cite{MaN} Let $(V,Y,{\bf 1})$ be a vertex algebra and let $v\in V$. We set $Y(v,z)=\sum_{n\in\mathbb{Z}}L(n)z^{-n-2}$. The vector $v$ is a {\em conformal vector} if 

\vspace{0.1cm} 

\noindent (i) there is a complex number $c$ such that $$[L(m),L(n)]=(m-n)L(m+n)+\frac{m^3-m}{12}\delta_{m+n,0}cId,$$

\vspace{0.1cm}

\noindent (ii) $L(-1)u=u_{-2}{\bf 1}$ for all $u\in V$. Here, $Y(u,z)=\sum_{m\in\mathbb{Z}}u_mz^{-m-1}$.

\vspace{0.1cm}

\noindent (iii) $L(0)$ is diagonalizable on $V$. 

\vspace{0.1cm}

\noindent The number $c$ is called the {\em central charge} of $v$. A vertex algebra endowed with a conformal vector $v$ is called a {\em conformal} vertex algebra of {\em rank }$c$. 
\end{dfn}
\begin{rmk} The condition (i) is equivalent to the conditions
$$L(n)v=\begin{cases}0~(n\geq 4),\\ \frac{c}{2}{\bf 1}~ (n=2),\\ 2v ~(n=0).\end{cases}$$
\end{rmk}

\vspace{0.2cm}

\noindent Recall that if $U$ is a vertex sub-algebra of $\overline{V_{B_{\lambda}}}$ that is generated by $\mathfrak{g}$ then $U$ is isomorphic to the rational vertex operator algebra $L_{\hat{g}}(1,0)$ associated to the Lie algebra $\mathfrak{g}$. Hence, we may assume that $L_{\hat{\mathfrak{g}}}(1,0)$ is a subvertex algebra of $\overline{V_{B_{\lambda}}}$. 

\vspace{0.2cm}

\noindent We now assume that 

\noindent i) $\omega$ is a conformal vector of $L_{\hat{g}}(1,0)$ such that $\omega\in \overline{V_{B_{\lambda}}}_{(2)}$, and 

\noindent ii) $c$ is a central charge of $L_{\hat{g}}(1,0)$. 

\vspace{0.2cm}

\noindent For convenience, we set $Y(\omega,z)=\sum_{m\in\mathbb{Z}}L(m)z^{-m-2}$. We let $h\in\mathfrak{h}$ and let $\{a_1,...,a_t\}$ be a linearly independent subset of $N_{\lambda}$ such that $L(1)h=0$, and $L(0)h=\alpha h$, $L(0)a_i=\beta_i a_i$, $h(0)a_i=(\lambda(h)-n_i)a_i$. Here, $\alpha,\beta_i\in \mathbb{Z}$, $n_i\geq 0$. Now, we set
\begin{eqnarray*}
\tilde{\omega}&&:=\omega+h(-2){\bf 1}+\sum_{i=1}^t\partial(a_i)(-2){\bf 1}\\
&&=\omega+h(-2){\bf 1}+2\sum_{i=1}^ta_i(-3){\bf 1}\in \overline{V_{B_{\lambda}}}_{(2)};\\
Y(\tilde{\omega},z)&&=\sum_{n\in\mathbb{Z}}\tilde{L}(n)z^{-n-2}.
\end{eqnarray*}
Hence, $\tilde{L}(n)=L(n)-(n+1)h(n)+n(n+1)\sum_{i=1}^ta_i(n-1)$. In particular, we have
\begin{eqnarray*}
\tilde{L}(-1)&&=L(-1),\\
\tilde{L}(0)&&=L(0)-h(0).
\end{eqnarray*}
Since $\tilde{\omega}\in \overline{V_{B_{\lambda}}}_{(2)}$, we then have that $\tilde{L}(n)\tilde{\omega}=0$ for all $n\geq 3$.

\vspace{0.2cm}

\noindent Now, we will find necessary and sufficient conditions that $\tilde{\omega}$ is a conformal vector of $\overline{V_{B_{\lambda}}}$.
\begin{lem} $\tilde{L}(2)\tilde{\omega}\in\mathbb{C}{\bf 1}$ if and only if $\beta_i=\lambda(h)-n_i$ for all $i\in\{1,...,t\}$. Moreover, if $\beta_i=\lambda(h)-n_i$ for all $i\in\{1,...,t\}$ then $\tilde{L}(2)\tilde{\omega}=\frac{c}{2}{\bf 1}-6\langle h,h\rangle{\bf 1}$.
\end{lem}
\begin{proof}
Since
\begin{eqnarray*}
\tilde{L}(2)\omega&&=(L(2)-3h(2)+6\sum_{i=1}^ta_i(1))\omega\\
&&=\frac{c}{2}{\bf 1}-3h(2)\omega+6\sum_{i=1}^ta_i(1)\omega\\
&&=\frac{c}{2}{\bf 1}+3(\sum_{j\geq 0}\frac{D^j}{j!}(-1)^j\omega(j+2)h)+6\sum_{i=1}^t\sum_{j\geq 0}\frac{D^j}{j!}(-1)^j\omega(1+j)a_i\\
&&=\frac{c}{2}{\bf 1}+6\sum_{i=1}^tL(0)a_i\\
&&=\frac{c}{2}{\bf 1}+6\sum_{i=1}^t\beta_i a_i,\\
\tilde{L}(2)h(-2){\bf 1}&&=-3h(2)h(-2){\bf 1}+6\sum_{i=1}^ta_i(1)h(-2){\bf 1}\\
&&=-6\langle h,h\rangle {\bf 1}-6\sum_{i=1}^t(\lambda(h)-n_i)a_i,\\
\end{eqnarray*}
and 
\begin{eqnarray*}
\tilde{L}(2)a_i(-3){\bf 1}&&=L(2)a_i(-3){\bf 1}-3h(2)a_i(-3){\bf 1}-6a_i(1)a_i(-3){\bf 1}\\
&&=\omega_3a_i(-3){\bf 1}-3(\lambda(h)-n_i)a_i\\
&&=\sum_{j\geq 0}{3\choose j}(\omega_ja)(-j){\bf 1}-3(\lambda(h)-n_i)a_i\\
&&=3L(0)a_i-3(\lambda(h)-n_i)a_i\\
&&=3(\beta_i-\lambda(h)+n_i)a_i,
\end{eqnarray*}
we have 
\begin{eqnarray*}
&&\tilde{L}(2)\tilde{\omega}\\
&&=\frac{c}{2}{\bf 1}+6\sum_{i=1}^t\beta_i a_i+(-6\langle h,h\rangle {\bf 1}-6\sum_{i=1}^t(\lambda(h)-n_i)a_i)+2\sum_{i=1}^t3(\beta_i-\lambda(h)+n_i)a_i\\
&&=\frac{c}{2}{\bf 1}-6\langle h,h\rangle{\bf 1}+12\sum_{i=1}^t(\beta_i-\lambda(h)+n_i)a_i.
\end{eqnarray*}
Clearly, $\tilde{L}(2)\tilde{\omega}\in\mathbb{C}{\bf 1}$ if and only if $\beta_i=\lambda(h)-n_i$ for all $i\in\{1,...,t\}$. In addition, if $\beta_i=\lambda(h)-n_i$ for all $i\in\{1,...,t\}$, we have $\tilde{L}(2)\tilde{\omega}=\frac{c}{2}{\bf 1}-6\langle h,h\rangle{\bf 1}$.
\end{proof}

\begin{lem} Assume that $\beta_i=\lambda(h)-n_i$ for all $i\in\{1,...,t\}$. Then $\tilde{L}(0)\tilde{\omega}=2\tilde{\omega}$ if and only if $\sum_{i=1}^ta_i(-3){\bf 1}=\frac{1}{2}\sum_{i=1}^t(L(-1)a_i)(-2){\bf 1}$.
\end{lem}
\begin{proof} Observe that
\begin{eqnarray*}
\tilde{L}(0)\omega&&=L(0)\omega-h(0)\omega\\
&&=2\omega-(-\omega(0)h+D\omega(1)h-\frac{D^2}{2!}\omega(2)h)\\
&&=2\omega-(-L(-1)h+\alpha D(h))\\
&&=2\omega-(\alpha-1)h(-2){\bf 1},\\
\tilde{L}(0)h(-2){\bf 1}&&=L(0)h(-2){\bf 1}-h(0)h(-2){\bf 1}\\
&&=\omega(1)h(-2){\bf 1}\\
&&=\sum_{j=0}^1{1\choose j}(\omega_jh)(-1-i){\bf 1}\\
&&=(L(-1)h)+(L(0)h)(-2){\bf 1}\\
&&=(1+\alpha)h(-2){\bf 1}.
\end{eqnarray*}
Also, by commutator formula and iterate formula, we have
\begin{eqnarray*}
\tilde{L}(0)a_i(-3){\bf 1}&&=L(0)a_i(-3){\bf 1}-h(0)a_i(-3){\bf 1}\\
&&=(\omega_0a_i)(-2){\bf 1}+(\omega_1a_i)(-3){\bf 1}-(h(0)a_i)(-3){\bf 1}\\
&&=\omega_0a_i(-2){\bf 1}+(\lambda(h)-n_i)a_i(-3){\bf 1}-(\lambda(h)-n_i)a_i(-3){\bf 1}\\
&&=(L(-1)a_i)(-2){\bf 1}.\\
\end{eqnarray*}
Therefore, we have
\begin{eqnarray*}
\tilde{L}(0)\tilde{\omega}&&=2\omega-(\alpha-1)h(-2){\bf 1}+(1+\alpha)h(-2){\bf 1}+2\sum_{i=1}^t(L(-1)a_i)(-2){\bf 1}\\
&&=2\omega+2h(-2){\bf 1}+2\sum_{i=1}^t(L(-1)a_i)(-2){\bf 1}.\\
\end{eqnarray*}
Hence, 
$\tilde{L}(0)\tilde{\omega}=2\tilde{\omega}$ if and only if $\sum_{i=1}^ta_i(-3){\bf 1}=\frac{1}{2}\sum_{i=1}^t(L(-1)a_i)(-2){\bf 1}$.
\end{proof}
\begin{rmk} If we assume that $L(-1)a=a(-2){\bf 1}$ for all $a\in N_{\lambda}$ then by iterate formula, we have $(L(-1)a)(-2){\bf 1}=2a(-3){\bf 1}$. In addition, we have $\tilde{L}(0)\tilde{\omega}=2\tilde{\omega}$.
\end{rmk}
\begin{thm} Let $(A,*)$ be a finite dimensional commutative associative algebra with the identity $\hat{1}$, and let $B$ be a finite dimensional vertex $A$-algebroid such that 

\vspace{0.1cm}

\noindent (i) $B$ is a simple Leibniz algebra that is not a Lie algebra. Consequently, there exists a simple Lie algebra $\mathfrak{g}$ such that $B=\mathfrak{g}+\partial(A)$ and $\partial(A)$ is a simple $\mathfrak{g}$-module;

\vspace{0.1cm}

\noindent (ii) there exist $u,v\in\mathfrak{g}$ such that $u_1v\neq 0$, $Ker(\partial)=\mathbb{C}\hat{1}$ and $N\cdot B\subset \partial(A)$. 

\vspace{0.1cm}

\noindent Hence, there exists a dominant integral $\lambda$ such that $B\cong B_{\lambda}$ as vertex algbroid. For simplicity, from now on we will assume that $$A:=A_{\lambda}=\mathbb{C}\hat{1}+N_{\lambda},~B:=B_{\lambda}=\mathfrak{g}+\partial(N_{\lambda}).$$ Here $\mathfrak{g}$ is a Lie algebra of $ADE$-type that is not $E_8$.  

\vspace{0.2cm}

\noindent We will denote the Killing form of $\mathfrak{g}$ by $\langle ~,~\rangle$. We let $\mathfrak{h}$ be a Cartan subalgebra of $\mathfrak{g}$ and identify $\mathfrak{h}$ with $\mathfrak{h}^*$ by means of the Killing form $\langle~,~\rangle$. Let $\Delta$ be the root system of $\mathfrak{g}$ with respect to $\mathfrak{h}$. We have the root space decomposition $\mathfrak{g}=\mathfrak{h}\oplus\sum_{\alpha\in\Delta}\mathfrak{g}_{\alpha}$. We let $\theta$ be the highest root. Hence $\langle\theta,\theta\rangle=2$. We will fix vectors $e_{\theta}\in\mathfrak{g}_{\theta}$, $f_{\theta}\in\mathfrak{g}_{-\theta}$ such that $(e_{\theta})_1f_{\theta}=\langle e_{\theta},f_{\theta}\rangle\hat{1}=\hat{1}$. 

\vspace{0.1cm}

\noindent Let $\omega$ be a conformal vector of $L_{\hat{g}}(1,0)$ such that $\omega\in \overline{V_{B_{\lambda}}}_{(2)}$, and set $Y(\omega,z)=\sum_{m\in\mathbb{Z}}L(m)z^{-m-2}$. Also, we let $h\in\mathfrak{h}$ and let $\{a_1,...,a_t\}$ be a basis of $N_{\lambda}$ such that $L(1)h=0$, 
$L(0)a_i=h(0)a_i=(\lambda(h)-n_i)a_i$, and $L(-1)a_i=a_i(-2){\bf 1}$ for all $i\in\{1,...,t\}$. Here, $n_i\geq 0$.

\vspace{0.1cm}

\noindent If we set $\tilde{\omega}=\omega+h(-2){\bf 1}+\sum_{i=1}^t\partial(a_i)(-1){\bf 1}$, then $\tilde{\omega}$ is a conformal vector of $\overline{V_{B_{\lambda}}}$. In addition, $\overline{V_{B_{\lambda}}}$ is a $C_2$-cofinite indecomposable non-simple conformal vertex algebra of rank $(c-12\langle h,h\rangle){\bf 1}$. 
\end{thm}
\begin{proof} To show that $\overline{V_{B_{\lambda}}}$ is a conformal vertex algebra of rank $(c-12\langle h,h\rangle){\bf 1}$, it is enough to show that $\tilde{L}(-1)u=u(-2){\bf 1}$ for all $u\in \overline{V_{B_{\lambda}}}$ and $\tilde{L}(0)$ is diagonalizable on $\overline{V_{B_{\lambda}}}$. Notice that 
\begin{eqnarray*}
&&\tilde{L}(-1)a'=L(-1)a'=a'(-2){\bf 1}=\partial(a')=D(a'),\\ 
&&\tilde{L}(-1)g=L(-1)g=g(-2){\bf 1}=D(g)
\end{eqnarray*} 
for all $a'\in N_{\lambda}$, $g\in \mathfrak{g}$. In addition, we have 
\begin{eqnarray*}
\tilde{L}(-1)a'(-n){\bf 1}&&=L(-1)a'(-n){\bf 1}\\
&&=\omega_0a'(-n){\bf 1}\\
&&=(\omega_0 a')(-n){\bf 1}\\ 
&&=(L(-1)a')(-n){\bf 1}\\ 
&&=(Da')(-n){\bf 1}\\
&&=D(a'(-n){\bf 1})~\text{ and }\\
\tilde{L}(-1)g(-n){\bf 1}&&=(Dg)(-n){\bf 1}=D(g(-n){\bf 1})
\end{eqnarray*} 
for all $a'\in N_{\lambda},g\in\mathfrak{g}, n\in\mathbb{Z}$. Using the fact that 
\begin{eqnarray*}
\tilde{L}(-1)b(-n)v&&=L(-1)b(-n)v\\
&&=\omega_0b(-n)v\\
&&=b(-n)\omega_0(v)+(\omega_0b)(-n)v\\
&&=b(-n)L(-1)v+(L(-1)b)(-n)v\\
&&=b(-n)L(-1)v+(Db)(-n)v
\end{eqnarray*} for all $b\in B_{\lambda},~v\in\overline{V_{B_{\lambda}}}$, we can conclude that $$\tilde{L}(-1)u=u(-2){\bf 1}=D(u)\text{ for all }u\in \overline{V_{B_{\lambda}}}.$$

\vspace{0.1cm}

\noindent Since $\tilde{L}(0)=L(0)-h(0)$, and $L(0)$ and $h(0)$ act semisimply on $\overline{V_{B_{\lambda}}}$, we can conclude immediately that $\tilde{L}(0)$ acts semisimply on $\overline{V_{B_{\lambda}}}$ as well. Moreover, $\overline{V_{B_{\lambda}}}$ is a $C_2$-cofinite indecomposable non-simple conformal vertex algebra of rank $(c-12\langle h,h\rangle){\bf 1}$ as desired. 
\end{proof}

\section{Appendices}
\subsection{Background on Leibniz Algebras}
\begin{dfn}(\cite{DMS}, \cite{FM})\ \ 

\vspace{0.2cm}

\noindent (i) A {\em left Leibniz algebra} $\mathfrak{L}$ is a $\mathbb{C}$-vector space equipped with a bilinear map $[~,~]:\mathfrak{L}\times\mathfrak{L}\rightarrow\mathfrak{L}$ satisfying the Leibniz identity $[a,[b,c]]=[[a,b],c]+[b,[a,c]]$ for all $a,b,c\in\mathfrak{L}$.

\vspace{0.2cm}

\noindent (ii) Let $\mathfrak{L}$ be a left Leibniz algebra over $\mathbb{C}$. Let $I$ be a subspace of $\mathfrak{L}$. $I$ is a {\em left} (respectively, {\em right}) {\em ideal} of $\mathfrak{L}$ if $[\mathfrak{L}, I]\subseteq I$ (respectively, $[I,\mathfrak{L}]\subseteq I$). $I$ is an {\em ideal} of $\mathfrak{L}$ if it is both a left and a right ideal. 
\end{dfn}

\begin{ex} We define $Leib(\mathfrak{L})=Span\{~[u,u]~|~u\in\mathfrak{L}~\}=Span\{[u,v]+[v,u]~|~u,v\in\mathfrak{L}\}$. $Leib(\mathfrak{L})$ is an ideal of $\mathfrak{L}$. Moreover, for $v,w\in Leib(\mathfrak{L})$, $[v,w]=0$. 
\end{ex}

\begin{dfn}\cite{DMS} Let $(\mathfrak{L}, [~,~])$ be a left Leibniz algebra. The series of ideals $$...\subseteq \mathfrak{L}^{(2)}\subseteq \mathfrak{L}^{(1)}\subseteq \mathfrak{L}$$ where $\mathfrak{L}^{(1)}=[\mathfrak{L},\mathfrak{L}]$, $\mathfrak{L}^{(i+1)}=[\mathfrak{L}^{(i)},\mathfrak{L}^{(i)}]$ is called the {\em derived series} of $\mathfrak{L}$.  A left Leibniz algebra $\mathfrak{L}$ is {\em solvable} if $\mathfrak{L}^{(m)}=0$ for some integer $m\geq 0$. As in the case of Lie algebras, any left Leibniz algebra $\mathfrak{L}$ contains a unique maximal solvable ideal $rad(\mathfrak{L})$ called the  the {\em radical} of $\mathfrak{L}$ which contains all solvable ideals.
\end{dfn}

\begin{ex} $Leib(\mathfrak{L})$ is a solvable ideal.
\end{ex}

\begin{dfn}\cite{DMS}\ \  

\vspace{0.2cm}

\noindent (i) A left Leibniz algebra $\mathfrak{L}$ is {\em simple}  if $[\mathfrak{L},\mathfrak{L}]\neq Leib(\mathfrak{L})$, and $\{0\}$, $Leib(\mathfrak{L})$, $\mathfrak{L}$ are the only ideals of $\mathfrak{L}$.

\vspace{0.2cm}

\noindent (ii) A left Leibniz algebra $\mathfrak{L}$ is said to be {\em semisimple} if $rad(\mathfrak{L})=Leib(\mathfrak{L})$. 
\end{dfn}
\begin{prop}(\cite{Ba2}, \cite{DMS}) Let $\mathfrak{L}$ be a left Leibniz algebra. 

\vspace{0.2cm}

\noindent (i) There exists a subalgebra $S$ which is a semisimple Lie algebra of $\mathfrak{L}$ such that $\mathfrak{L}=S \dot{+} rad(\mathfrak{L})$. As in the case of a Lie algebra, we call $S$ a Levi subalgebra or a Levi factor of $\mathfrak{L}$.

\vspace{0.2cm}

\noindent (ii) If $\mathfrak{L}$ is a semisimple Leibniz algebra then $\mathfrak{L}=(S_1\oplus S_2\oplus...\oplus S_k)\dot{+}Leib(\mathfrak{L})$, where $S_j$ is a simple Lie algebra for all $1\leq j\leq k$. Moreover, $[\mathfrak{L},\mathfrak{L}]=\mathfrak{L}$. 

\vspace{0.2cm}

\noindent (iii) If $\mathfrak{L}$ is a simple Leibniz algebra, then there exists a simple Lie algebra $S$ such that $Leib(\mathfrak{L})$ is an irreducible module over $S$ and $\mathfrak{L}=S\dot{+}Leib(\mathfrak{L})$.
\end{prop}

\begin{dfn} Let $\mathfrak{L}$ be a left Leibniz algebra. A left $\mathfrak{L}$-module is a vector space $M$ equipped with a $\mathbb{C}$-bilinear map $\mathfrak{L}\times M\rightarrow M; (u,m)\mapsto u\cdot m$ such that $([u,v])\cdot m=u\cdot (v\cdot m)-v\cdot(u\cdot m)$ for all $u,v\in \mathfrak{L},m\in M$. 

\vspace{0.2cm}

\noindent The usual definitions of the notions of submodule, irreducibility, complete reducibility, homomorphism, isomorphism, etc., hold for left Leibniz modules.
\end{dfn}

\begin{rmk} For any module $M$ over a Leibniz algebra $\mathfrak{L}$, $Leib(\mathfrak{L})$ acts as zero on $M$. 
\end{rmk}
\subsection{Background on Vertex Algebras} 

\ \

\begin{dfn}(\cite{B1}, \cite{FLM2}, \cite{LLi}) A {\em vertex algebra} is a vector space $V$ equipped with a linear map 
\begin{eqnarray*}
Y:V&&\rightarrow \End(V)[[x,x^{-1}]]\\
v&&\mapsto Y(v,x)=\sum_{n\in\mathbb{Z}}v_nx^{-n-1}\text{ where } v_n\in\End(V)
\end{eqnarray*} 
and equipped with a distinguished vector ${\bf 1}$, the {\em{vacuum vector}}, such that for $u,v\in V$, 
\begin{eqnarray*}
&&u_nv=0\text{ for $n$ sufficiently large},\\
&&Y({\bf 1},x)=1,\\
&&Y(v,x){\bf 1}\in V[[x]],\text{ and }\lim_{x\rightarrow 0}Y(v,x){\bf 1}=v
\end{eqnarray*}
and such that 
\begin{eqnarray*}
&&x_0^{-1}\delta\left(\frac{x_1-x_2}{x_0}\right)Y(u,x_1)Y(v,x_2)-x_0^{-1}\delta\left(\frac{x_2-x_1}{-x_0}\right)Y(v,x_2)Y(u,x_1)\\
&&=x_2^{-1}\delta\left(\frac{x_1-x_0}{x_2}\right)Y(Y(u,x_0)v,x_2)
\end{eqnarray*}
the {\em Jacobi identity}.
\end{dfn}

\vspace{0.2cm}

\noindent From the Jacobi identity we have Borcherds' commutator formula and iterate formula:
\begin{eqnarray}
&&{[u_m,v_n]}=\sum_{i\geq 0}{m\choose i}(u_iv)_{m+n-i},\\
&&(u_mv)_nw=\sum_{i\geq 0}(-1)^i{m\choose i}(u_{m-i}v_{n+i}w-(-1)^mv_{m+n-i}u_iw)
\end{eqnarray}
for $u,v,w\in V$, $m,n\in \mathbb{Z}$. 

\vspace{0.2cm}

\noindent We define a linear operator $D$ on $V$ by $D(v)=v_{-2}{\bf 1}$ for $v\in V$. Then 
\begin{eqnarray*}
&&Y(v,x){\bf 1}=e^{xD}v\text{ for }v\in V,\text{ and }\\
&&[D,v_n]=(Dv)_n=-nv_{n-1}\text{ for }v\in V, ~n\in\mathbb{Z}.
\end{eqnarray*} Moreover, for $u,v\in V$, we have $Y(u,x)v =e^{xD}Y(v,-x)u$ (skew-symmetry).

\vspace{0.2cm} 

\noindent A vertex algebra $V$ equipped with a $\mathbb{Z}$-grading $V=\oplus_{n\in \mathbb{Z}}V_{(n)}$ is called a {\em $\mathbb{Z}$-graded vertex algebra} if ${\bf 1}\in V_{(0)}$ and if $u\in V_{(k)}$ with $k\in\mathbb{Z}$ and for $m,n\in\mathbb{Z}$, $u_mV_{(n)}\subseteq V_{(k+n-m-1)}$. 

\vspace{0.2cm}

\noindent A $\mathbb{N}$-graded vertex algebra is defined in the obvious way.

\begin{prop}\cite{GMS} Let $V=\oplus_{n\in \mathbb{N}}V_{(n)}$ be an $\mathbb{N}$-graded vertex algebra. Then 

\vspace{0.2cm}

\noindent (i) $V_{(0)}$ is a commutative associative algebra with ${\bf 1}$ as identity where $a* a'=a_{-1}a'$ for $a, a'\in V_{(0)}$.

\noindent (ii) $V_{(1)}$ is a Leibniz algebra with $[b,b']=b_0b'$ for $b,b'\in V_{(1)}$.

\noindent (iii) The graded vector space $V_{(0)}\oplus V_{(1)}$ is a 1-truncated conformal algebra with the linear map $D:V_{(0)}\rightarrow V_{(1)}$, and bilinear operations $(u,v)\mapsto u_iv$ for $i=0,1$ on $V_{(0)}\oplus V_{(1)}$.

\noindent (iv) Moreover, $V_{(1)}$ is a vertex $V_{(0)}$-algebroid with $a\cdot b=a_{-1}b$ for $a\in V_{(0)}, b\in V_{(1)}$.
\end{prop}

\begin{dfn}(\cite{FLM2}, \cite{LLi}) An ideal of a vertex algebra $V$ is a subspace $I$ such that $u_nw\in I$ and $w_nu\in I$ for $u\in V$, $w\in I$ and $n\in\mathbb{Z}$. 
\end{dfn} 

\noindent Notice that for any ideal $I$ of $V$, we have $D(w)=w_{-2}{\bf 1}\in I$. Hence, under the condition that $D(I)\subseteq I$, the left ideal condition $v_nw\in I$ for all $v\in V$, $w\in I$, $n\in\mathbb{Z}$ is equivalent to the right ideal condition $w_mv\in I$ for all $v\in V$, $w\in I$, $m\in\mathbb{Z}$.

\vspace{0.2cm}

\noindent For a subset $S$ of a vertex algebra $V$, we denote by $(S)$ the smallest ideal of $V$ containing $S$. It was shown in Corollary 4.5.10 of \cite{LLi} that 
$$(S)=Span\{~v_nD^i(u)~|~v\in V, n\in\mathbb{Z}, i\geq 0, u\in S\}.$$

\vspace{0.2cm}

\begin{dfn}  For a vertex algebra $V$, we define $C_2(V)=Span\{u_{-2}v~|~u,v\in V\}$. The vertex algebra $V$ is said to satisfy {\em  the $C_2$-condition} if $V/C_2(V)$ is finite dimensional
\end{dfn}

\begin{prop}\cite{Z} For $u,v\in V$, $n\geq 2$, $D(v)\in C_2(V)$ and $u_{-n}v\in C_2(V)$.
\end{prop}

\begin{dfn}\cite{LLi} A $V$-{\em module} is a vector space $W$ equipped with a linear map $Y_W$ from $V$ to $(\End W)[[x,x^{-1}]]$ where $Y_W(v,x)=\sum_{n\in\mathbb{Z}}v_nx^{-n-1}$ for $v\in V$ such that for $u,v\in V$, $w\in W$, 
\begin{eqnarray*}
&&u_nw=0\text{ for $n$ sufficiently large},\\
&&Y_W({\bf 1},x)=1,\\
&&x_0^{-1}\delta\left(\frac{x_1-x_2}{x_0}\right)Y_W(u,x_1)Y_W(v,x_2)-x_0^{-1}\delta\left(\frac{x_2-x_1}{-x_0}\right)Y_W(v,x_2)Y_W(u,x_1)\\
&&=x_2^{-1}\delta\left(\frac{x_1-x_0}{x_2}\right)Y_W(Y(u,x_0)v,x_2).
\end{eqnarray*}
\end{dfn}
\begin{dfn} Let $V=\oplus_{n=0}^{\infty}V_{(n)}$ be a $\mathbb{N}$-graded vertex algebra. A $\mathbb{N}$-graded $V$-module is a $V$-module $M$ equipped with a $\mathbb{N}$-grading $M=\oplus_{n=0}^{\infty}M_{(n)}$ such that $v_mM_{(n)}\subset M_{(n+p-m-1)}$ for $v\in V_{(p)}$, $p,n\in\mathbb{N}$, $m\in\mathbb{Z}$. For convenient, for a nonzero $\mathbb{N}$-graded $V$-module $M=\oplus_{n=0}^{\infty}M_{(n)}$, we may 
assume that $M_{(0)}\neq \{0\}$.
\end{dfn}

\begin{prop}\label{nil}\cite{LLi} Let $V$ be a vertex algebra.

\noindent (i) For $v\in V$, $v$ is weakly nilpotent if and only if $(v_{-1})^r{\bf 1}=0$ for some $r>0$. 

\noindent (ii) If $u\in V$ such that $u_nu=0$ for all $n\geq 0$, then $Y((u_{-1})^r{\bf 1},z)=Y(u,z)^r$ for $r>0$.

\noindent (iii) Let $(W,Y_W)$ be a $V$-module. Let $u\in V$ such that $u_nu=0$ for all $n\geq 0$, we have $Y_W((u_{-1})^r{\bf 1},z)=Y_W(u,z)^r$ for $r>0$. 
\end{prop}


\end{document}